\documentclass[12pt,a4paper]{article}
\usepackage[british]{babel}
\usepackage[a4paper,top=2cm,bottom=2cm,left=2.5cm,right=2.5cm,marginparwidth=1.75cm]{geometry}

\usepackage{xcolor} 
\definecolor{changes}{RGB}{220,50,47} 

\usepackage[style=ieee, backend=biber]{biblatex} 
\usepackage{amsmath}
\usepackage{amssymb}
\usepackage{amsfonts}
\usepackage{mathrsfs}
\usepackage{graphicx}
\usepackage[colorlinks=true, allcolors=blue]{hyperref}
\usepackage{orcidlink}
\usepackage[title]{appendix}

\usepackage{booktabs} 
\usepackage{threeparttable} 
\usepackage{algorithm}
\usepackage{algorithmicx}
\usepackage{algpseudocode}
\usepackage{listings}
\usepackage{enumitem}
\usepackage{chngcntr}
\usepackage{lipsum}
\usepackage{authblk}
\usepackage[T1]{fontenc}    
\usepackage{csquotes}       
\usepackage{diagbox}
\usepackage{comment}
\usepackage{siunitx}
\usepackage{helvet}  
\usepackage{lmodern}  
\usepackage{multirow}
\usepackage{placeins}
\usepackage{adjustbox}

\usepackage[justification=raggedright, singlelinecheck=false]{caption}
\usepackage{subcaption}
\usepackage{float}   

\newcommand{\apanum}[1]{\textbf{#1}} 
\DeclareFieldFormat[article]{volume}{\apanum{#1}} 

\usepackage{setspace}
\usepackage{titlesec}
\titleformat{\section} 
{\normalfont\Large\bfseries}{\thesection.}{1em}{}



\title{A Scalable Approach for Transient Thermal Modeling of Automotive Power Electronics}

\author{Neelakantan Padmanabhan}
\affil{\small ZF Active Safety and Electronics US Inc, \\ Livonia, MI, USA}
\affil{\small neel.padmanabhan@zf.com}

\begin{document}
	\maketitle
\footnotetext{\textcopyright 2025 SAE International. This work has been published in SAE Technical papers as: Padmanabhan, N., "A Scalable Approach for Transient Thermal Modeling of Automotive Power Electronics," SAE Technical Paper Series, January 1, 2025, https://doi.org/10.4271/2025-01-5073. This arXiv version corresponds to the author accepted manuscript. The final version of record is available at https://doi.org/10.4271/2025-01-5073}
\begin{abstract}
Efficient thermal management is critical for the reliability and performance of power electronics systems in automotive applications. This work presents a computationally efficient modeling approach for transient thermal simulation of power electronic systems, with a focus on inverter modules using multiple MOSFETs mounted on a printed circuit board assembly (PCBA). A case study of an inverter module comprising six MOSFETs arranged as high-side and low-side pairs for a three phases system mounted on a PCBA, attached to a heat sink is considered. Computational fluid dynamic (CFD) simulations in Ansys\textsuperscript{\textregistered} Icepak\textsuperscript{\texttrademark} are performed considering different heat transfer mechanisms, including natural convection, forced convection at constant velocity, and forced convection with varying flow velocity. A transient thermal model is developed using the Lumped Parameter Linear Superposition (LPLSP) method, a hybrid approach that combines lumped parameter modeling with the principle of linear superposition to capture transient thermal behavior efficiently. Temperatures of the components from the simulations are compared with temperatures from the LPLSP model and temperatures from a Linear Time Invariant (LTI) based reduced order model (ROM) developed for this system. It is observed that the LPLSP model is able to model a wide range of use cases very accurately with error of less than $5 \%$. This method enables rapid thermal performance evaluation of power electronics systems that have very fast transients in component level power dissipation and variations in ambient conditions, making it particularly well-suited for early-stage design iterations and long-duration mission profile simulations. The approach offers a practical path to reducing development cycles for automotive power electronics design. 
\end{abstract}


\section{Introduction}
Power electronics systems which include devices such as inverters, rectifiers, and power supplies, are fundamental in all modern automotive applications, including electric vehicles, internal combustion (IC) engine controls, and vehicle automation systems. Thermal management of power electronics is a critical challenge, as excessive heat can degrade performance, reduce the lifespan of components, and lead to system failure. Accurate thermal simulation is vital in the design and optimization of these systems \cite{neel_ecce}. However, conventional thermal simulations via computational fluid dynamics (CFD), are computationally expensive and time-consuming, especially for complex, large-scale systems \cite{neel_ecce,neel_cfd2}. Power electronics systems comprise of several components including semiconductors, capacitors, shunts, inductors, and diodes which generate heat through various mechanisms including conduction losses, switching losses, diode losses. 
Thermal simulations must account for the complex interactions between electrical, thermal, and fluid dynamics, as well as the detailed geometrical and material properties of the system. This complexity makes accurate simulations computationally expensive. Therefore, reduced order models (ROM) have emerged as an attractive alternative, offering a balance between accuracy and efficiency. Various ROM techniques have been developed to reduce the system’s degrees of freedom while preserving its essential thermal characteristics. Some of these approaches include Krylov-subspace models, data driven models including proper orthogonal decomposition (POD) and machine learning models as well as compact thermal networks. \\
The Krylov subspace approach employs an iterative process to approximate solutions of large systems by projecting the system onto a lower dimensional subspace. Krylov subspace ROMs are effective in systems whose thermal behavior can be approximated by linear models \cite{krylov1,krylov2}, such as in linear time invariant (LTI) systems. Krylov-subspace ROMs are available in commercial CFD tools like Ansys\textsuperscript{\textregistered} Icepak\textsuperscript{\texttrademark}. However, there is limited research on application of the method for systems with nonlinearities and coupled thermal-electrical behaviors. The accuracy of the ROMs depend heavily on choice of expansion points which requires trial-and-error \cite{krylov3,krylov4}. Further, construction and maintenance of orthonormal basis becomes computationally intensive as the subspace grows \cite{krylov5}. A similar projection based reduction approach is Balanced Truncation that approximates a system by retaining the most balanced modes, where the system’s energy is concentrated. This method is known for providing accurate approximations in LTI systems, but it is prohibitively expensive for very large scale systems and performs poorly in non-linear systems \cite{BT1} . Additionally, the stability of the ROMs for second-order systems are not guaranteed \cite{BT2}. \\
Proper Orthogonal Decomposition (POD) is one of the most widely used ROM methods for thermal simulations. Berkooz et al. \cite{pod_theory} provide foundational theory on POD for turbulence, which supports its adoption in thermal simulations. POD identifies the most important modes of the system by analyzing its dynamic response and then represents the system in a reduced subspace \cite{neel_pod, neel_chapter}. Zhu et al. \cite{zhu2020proper1, zhu2020proper2} present a POD framework with explicit and implicit time-marching schemes for transient nonlinear heat conduction, demonstrating high accuracy and notable reductions in computational cost. Further applications include system-level thermal modeling of shipboard power electronics cabinets using POD-derived compact models from CFD data, and similar strategies for data center thermal management \cite{haider2008proper}. However, despite its effectiveness in systems with dominant thermal modes, POD may offer limited computational gains in scenarios with strong nonlinearity or time-varying behavior, where the system's energy is more broadly distributed across many modes. \\
Over the past decade machine learning has substantially advanced thermal modeling in power electronics by enabling real-time prediction of faults, efficiency and remaining useful life. However, these data-driven approaches are highly dependent on quality and quantity of data and suffer from issues like limited or noisy datasets, overfitting and ability to generalize. Further, it is a black box approach that lacks interpretability, guarantees of stability and general robustness, which is a safety concern for critical power electronics systems \cite{ML1,ML2,ML3}. \\
Compact thermal network (CTN) models represent thermal systems as network of resistance and capacitance (RC), capturing the flow of heat through various components. This approach typically employs Foster and Cauer RC circuits to represent transient thermal impedance of devices \cite{fostercauer}. Foster models, derived by fitting thermal impedance curves, offer simplicity and ease of parameter extraction but lack physical correspondence to device layers. In contrast, Cauer models are constructed as laddered RC elements and map directly to physical thermal paths, enabling improved integration with heatsinks or PCBs \cite{CTN1,CTN2,CTN3,CTN4}. Reviews in Insulated Gate Bipolar Transistor (IGBT) modules highlight the necessity of both Foster and Cauer types, and stress that many conventional studies omit thermal coupling between chips. Common gaps include under representation of multi source thermal coupling, reliance on extensive measurements for parameter extraction and assumptions of linear thermal properties despite known temperature dependencies. These limitations constrain accurate and scalable thermal simulation of high‑power devices, suggesting a need for semi‑empirical, nonlinear, and automation friendly thermal network models with broader experimental validation. \\
This paper builds upon the work presented by Padmanabhan \cite{neel_ieee}, which introduced a semi-empirical thermal model for transient power electronics systems, by applying the framework to a real-time use case involving an automotive power electronic circuit. The formulation is further extended to incorporate convective heat transfer with variable flow velocity, and the model's performance is compared with that of an impulse response based reduced order model for LTI systems. While most ROM approaches (e.g., Krylov-subspace, lumped parameter modeling) offer strong performance in linear or quasi-linear regimes, they tend to struggle with systems exhibiting time-varying boundary conditions and nonlinear behavior. POD methods, while effective in systems dominated by a few thermal modes, offer limited gains when the thermal energy is distributed across multiple spatial and temporal scales, as is common in automotive applications. In contrast, the LPLSP model proposed in this work offers an interpretable framework that maintains high accuracy even in the presence of nonlinear and time-dependent heat transfer mechanisms. Unlike traditional CTN models, it captures thermal interactions between various sources and sinks without relying on exhaustive measurements. Moreover, it circumvents the stability and interpretability challenges of machine learning models, while providing computational speeds that are orders of magnitude faster than full CFD simulations. A sensitivity analysis is also conducted to evaluate the effects of physical parameters of the system on the modeled temperature. Additionally, the applicability of the modeling approaches to a range of heat transfer scenarios is examined.
\section{Case Study} \label{casestudy}
\begin{figure}
	\centering
	\includegraphics[width=\linewidth]{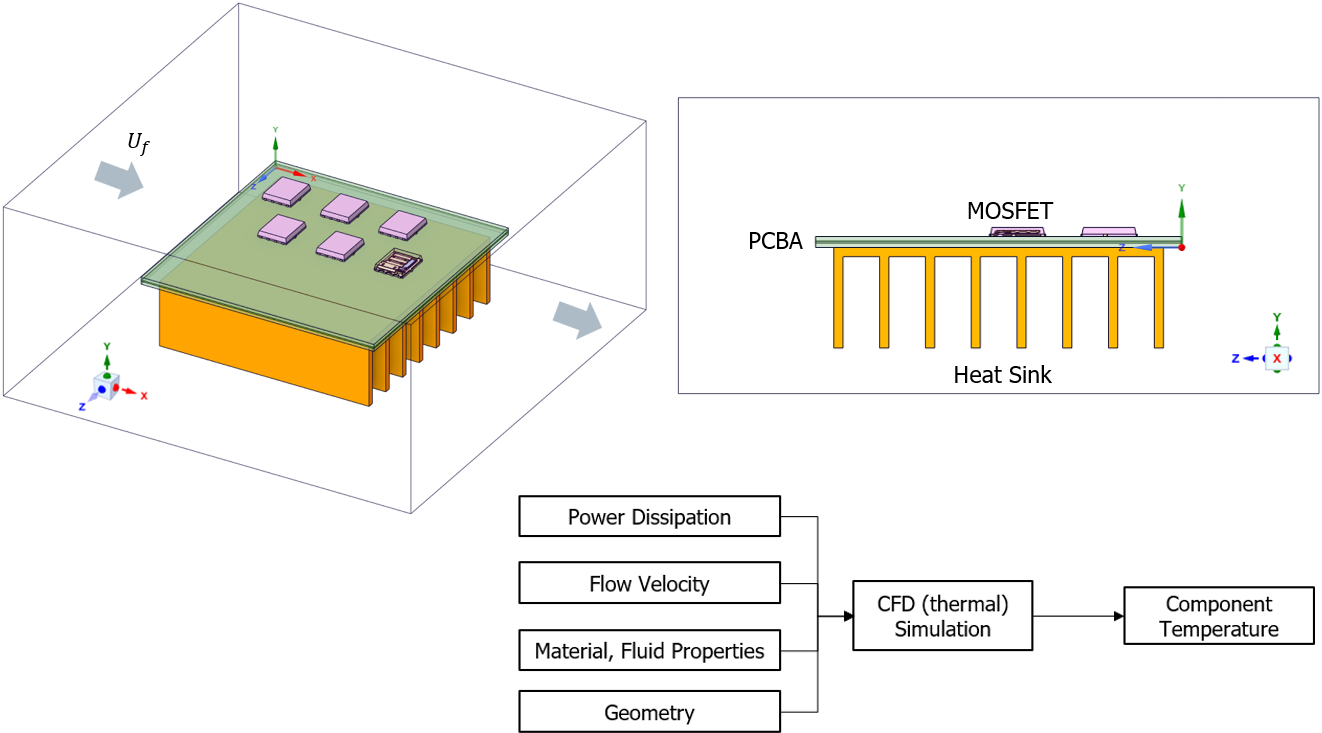}
	\caption{An illustration of the CFD simulation setup consisting of the MOSFETs, PCBA, and heat sink. Open boundary conditions are applied in all directions for natural convection and inlet boundary conditions are applied in $x$ direction for forced convection.}
	\label{fig:T_comp}
\end{figure}

The case study involves an inverter module comprising six MOSFETs mounted on a printed circuit board assembly (PCBA). An illustration of the simulation configuration presented in Fig.~\ref{fig:T_comp}. Heat is dissipated at the MOSFET junction due to conduction losses, switching losses, gate losses and diode losses. The heat is transferred from the MOSFET package to the PCBA and subsequently dissipated to the environment with a finned heat sink. For simplicity, the heat sink assumed to be in direct contact with the PCBA. This configuration is representative of a broad range of electronic circuits, making the approach described in this work, applicable to diverse electronic thermal management problems. The simulations are created in Ansys\textsuperscript{\textregistered} Icepak\textsuperscript{\texttrademark} under both natural convection and forced convection modes. Detailed mechanical CAD models, electrical CAD models and accurate material properties (Table~\ref{tab:material}) are used for the PCBA, MOSFETs, and heatsink. \\ %
Incompressible form of the Navier-Stokes equations along with energy, zero order turbulence equation and equation of state are solved by the solver (finite volume based). First and second order upwind schemes are applied for spatial discretization and fully implicit time integration scheme is used for temporal integration. For natural convection problems, the Boussinesq approximation is applied, while ideal gas equation of state is applied for the forced convection problem. Initial conditions are set to temperature of $\SI{20}{\degree} C$, pressure of $1 atm$ and standard earth gravity. The solver's default convergence criteria for the flow and energy equations are adopted. Open boundary conditions are specified for natural convection and an inlet $x-$velocity $U_f$ is specified for forced convection. The initial Rayleigh number is $3.1 \times 10^{5}$ and Prandtl number is $0.7$ for natural convection case. \\
Grid resolution is chosen such that the flow statistics become independent of the mesh density. Grid stretching and clustering are applied to concentrate mesh nodes in critical geometries, while non-critical regions have coarser meshing. In Icepak \textsuperscript{\texttrademark}, this is achieved through non-conformal meshing that enables separate meshing of assemblies. Multi-level meshing and slacks are applied to avoid mesh leakage and ensure smooth transitions between mesh densities. Meshing priorities are assigned to objects within assemblies and per object meshing is applied specifically to PCBA. Additionally, cut-cells are used for all objects. A mesh independence study was conducted to determine the final mesh size, with total grid points ranging from 0.3 million to 3 million.

\begin{table}[h] 
	\begin{center}
\begin{tabular}{|cccc|}
	\hline
	\multicolumn{1}{|c|}{Material}   & \multicolumn{1}{c|}{$\rho [kg/m^3]$} & \multicolumn{1}{c|}{$c_p [J/kg K]$} & $k [W/mK]$ \\ \hline
	\multicolumn{4}{|c|}{MOSFET}                                                                                               \\ \hline
	\multicolumn{1}{|c|}{Silicon}    & \multicolumn{1}{c|}{2300}            & \multicolumn{1}{c|}{793}            & 100        \\ \hline
	\multicolumn{1}{|c|}{Die Attach} & \multicolumn{1}{c|}{8500}            & \multicolumn{1}{c|}{50}             & 150        \\ \hline
	\multicolumn{1}{|c|}{Mold}       & \multicolumn{1}{c|}{1900}            & \multicolumn{1}{c|}{795}            & 0.8        \\ \hline
	\multicolumn{4}{|c|}{PCB}                                                                                                  \\ \hline
	\multicolumn{1}{|c|}{FR4}        & \multicolumn{1}{c|}{1250}            & \multicolumn{1}{c|}{1300}           & 0.35       \\ \hline
	\multicolumn{1}{|c|}{Copper}     & \multicolumn{1}{c|}{8933}            & \multicolumn{1}{c|}{387}            & 385        \\ \hline
	\multicolumn{4}{|c|}{Heat Sink}                                                                                            \\ \hline
	\multicolumn{1}{|c|}{Aluminum}   & \multicolumn{1}{c|}{2800}            & \multicolumn{1}{c|}{900}            & 205        \\ \hline
\end{tabular}
		\caption{Material properties of components used in the CFD simulation. $\rho$ represents the material density, $c_p$ the specific heat and $k$ the thermal conductivity.}
\label{tab:material}		
	\end{center}
\end{table}
\section{Lumped Parameter Linear Superposition Model}
The LPLSP method is empirically derived but grounded in the theory of heat transfer and can be related to the generalized form of the heat equation. The formulation process began with CFD simulations developed for specific individual cases. A specific form of the governing equations, as presented in the next section, was then fitted to the temperature data obtained from these simulations. The resulting model constants were optimized and subsequently analyzed to establish their physical relevance and consistency. This led to a more concrete form of the equations, applicable across the various simulated cases. The formulation is developed in steps: it is first presented for a single insulated body, then extended to conduction between two bodies, then extended to convection between one body and fluid, and finally generalized to include mixed modes of heat transfer, including convection under both constant and variable flow velocities. The transient form of heat equation considering constant properties at atmospheric pressures can be expressed as,
\begin{equation} \label{eqmain1}
    \frac{\partial T}{\partial t} = \alpha \nabla^2 T + \frac{P}{\rho c_p},
\end{equation}
where, $P$ is the heat source, $\rho, c_p ,k$ represent the density, specific heat capacity and thermal conductivity of the material, respectively, $T$ the temperature, $x_i,t$ the spatial and temporal terms, $\alpha= k/\rho c_p$ is the thermal diffusivity and $\nabla^2$ the Laplace operator. The transient temperature response can be expressed by integrating Eq. \eqref{eqmain1} with respect to time, 
\begin{equation} \label{eqmain}
    T(t)=D(t) + \frac{P}{\rho c_p} t + c_{int}. 
\end{equation}
where, $D(t)=\alpha \int \nabla^2 T dt$ represents the time integrated diffusion term and $c_{int}$ the integration constant which can be determined from the initial condition. This equation forms the basis of the proposed model. The model derivation involved fitting a specific form of this equation to the temperature data obtained from CFD simulations. The resulting model constants were optimized and then analyzed to relate them back to the physical parameters of the system.
\subsection{Model for insulated system}
For a single insulated cuboidal body, the transient temperature response is observed to be linear. This behavior can be derived from Eq.~\ref{eqmain} under the assumption that the temperature diffusion across a body is instantaneous, or equivalently, the spatial gradient of temperature is negligible. Under this assumption, the diffusion term $D(t)$ can be considered constant or zero, resulting in a simplified expression for temperature evolution as,
\begin{equation}
    T(t)=T^0 + T^L(t) = T^0 + \Bigl(\frac{P}{C} \Bigr) t,
\end{equation}
where, $T^0$ is the initial temperature, $P$ the total power applied to the body, $C=\rho c_p V$ the thermal capacitance of the body, and $V$ the volume of the body. Simulation results for the insulated body confirm that the temperature increases linearly over time, with a slope given by $T^L = (P/C)$.  
\subsection{Model for conduction}
\begin{figure}
    \centering
    \includegraphics[width=\linewidth]{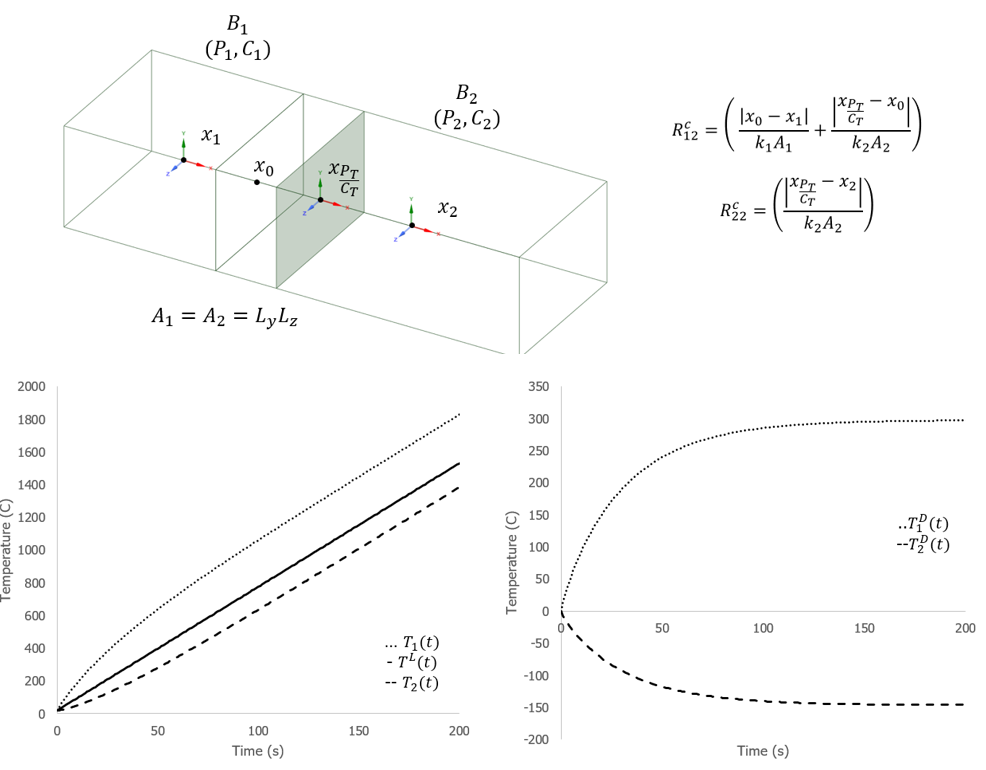}
    \caption{Conduction between two bodies insulated from the surrounding. Temperatures of the bodies $T_1(t), T_2(t)$, linear temperature $T^L(t)$ and deviation terms $T^D_1(t), T^D_2(t)$.}
    \label{fig:2Body}
\end{figure}

\begin{figure}
    \centering
    \includegraphics[width=\linewidth]{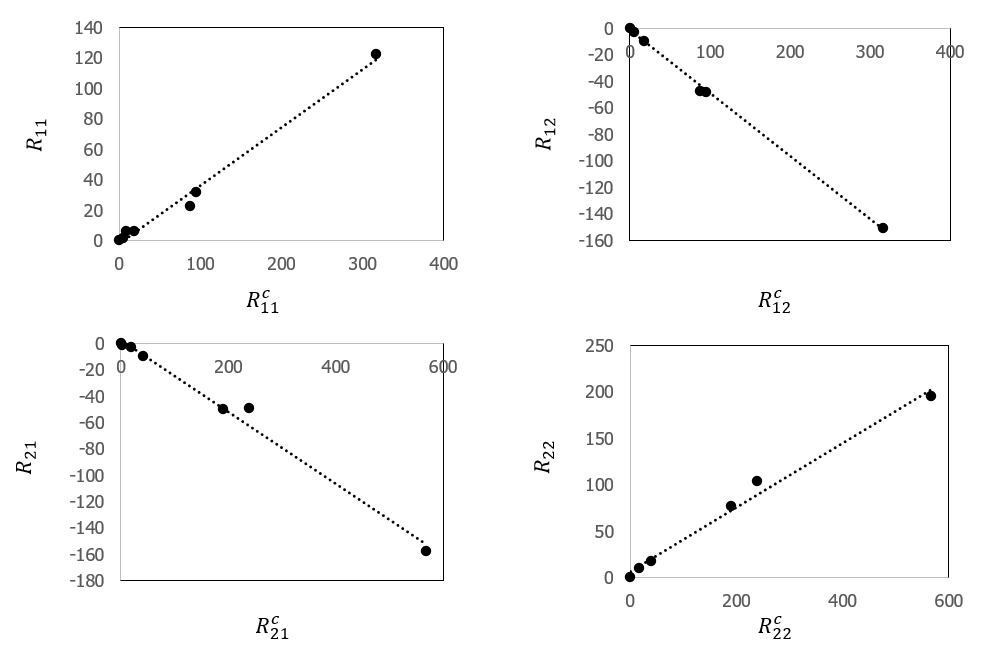}
    \caption{Thermal resistance estimated empirically and calculated analytically}
    \label{fig:Rcalc}
\end{figure}

Extending from a single insulated body to two cuboidal bodies insulated from the surroundings, it is observed that the spatial gradient can no longer be neglected due to heat conduction from source to sink body. Although the geometry is three-dimensional, the temperature gradient is observed to be significant only along the direction of heat flow. The temperatures at the centers of each body $T_i (t)$ are observed to be,
\begin{equation}
\begin{aligned}
    T _{1} (t)=T^0+T^L(t)+T^{D} _{1}(t),\\ 
    T _{2} (t)=T^0+T^L(t)+T^{D} _{2}(t).
\end{aligned}
\end{equation}
The linear temperature term is defined as $T^L(t)=(P_T / C_T) t$, where the ratio $P_T/C_T = \sum_{i=1} ^2 P_i / \sum_{i=1} ^2 C_i$ represents the total power input applied to both bodies divided by their combined thermal capacitance. The linear temperature curve $T^L(t)$ lies between temperatures $T_1(t)$ and $T_2(t)$, and can be computed a priori using the material properties, geometrical details of the bodies and input power dissipation. The spatial location (plane in this case) denoted by $x_{P_{T}/C_{T}}$ where $T^L(t)$ is observed, depends on the applied power dissipation, and thermal capacitance of each body. For bodies with similar thermal capacitances and power dissipation, this location typically lies near the midpoint between the two bodies. However, for bodies with differing capacitances, the linear temperature plane shifts closer to the body with the larger thermal capacitance. The term $T^{D}_i (t)$, models the deviation in temperature between the body $T^{s}_i(t)$, measured from the simulation and $T^L(t)$. This can be expressed as,
\begin{equation}
	T^{D}_i (t)= T^{s}_i (t) - T^L(t). 
\end{equation}
This deviation exhibits the characteristic transient response of an RC circuit \cite{RC_tr1,RC_tr2} and can be expressed as,
\begin{equation}
    T^{D} _{i}(t)=E_i (1-e^{-\theta_i t}).
\end{equation}
An illustration of this is presented in Fig.~\ref{fig:2Body}. The term $E_i$ is the magnitude of deviation between the stationary temperature of body $i$ and the linear temperature curve $T^L(t)$, and is defined as,
\begin{equation}
	E_i = |T^{s} _{i}(t_{s})-T^L(t_{s}))|, 
\end{equation}
where $T^{s} _i(t_s)$ is the simulated temperature of the body $i$ and $t_s$ is the time after which the temperature reaches a quasi-stationary state. The sign of $E_i$ is positive if $T^{s} _{i} > T^L$, and negative otherwise, indicating whether the body is hotter or cooler than the linear temperature curve. This term is modeled as,
\begin{equation}
    E_i = P_{i} R _{ij},
\end{equation}
where, $P_{i}$ is the applied power at body $i$ and $R _{ij}$ the thermal resistance between the location $x_i$ (where the temperature $T^{s}_i$ is measured) and the location $x_{P_{T}/C_{T}}$. The term $\theta_i$ represents the time constant which characterizes the rate of change of temperature. In systems with two bodies where each body may act as the source, sink or both, the deviation terms must consider the effects of each source on the other. Thus, the principle of superposition is applied and the deviation terms are expressed as,
\begin{equation}
    \begin{aligned}
        T^{D} _{1}(t) = (P_{1} R_{11} + P_{2} R_{12}) [1-e^{-\theta_1 t}], \\
        T^{D} _{2}(t) = (P_{1} R_{21} + P_{2} R_{22}) [1-e^{-\theta_2 t}].
    \end{aligned}
\end{equation}
The term $R_{11}$ represents the thermal resistance between the body $B_1(x_1)$ and $x_{P_{T}/C_{T}}$, when power $P_1 = 1$ is applied to it and power $P_2=0$ is applied to $B_2$. This term can be obtained by running a simulation where $P_1=1, P_2=0$ is applied to bodies $B_1, B_2$, and computed as follows,  
\begin{equation}
    R_{11} = \frac{T_1 ^{s} (t_s) - T^L (t_s)}{P_1}.
\end{equation}
$T_1 ^{s} (t_s)$ in this equation is the temperature of $B_1$ measured from the simulation. Similarly, 
\begin{equation}
\begin{aligned}
    R_{12} = (T_1 ^{s} (t_s) - T^L (t_s))/P_2, \\
    R_{21} = (T_2 ^{s} (t_s) - T^L (t_s))/P_1, \\
    R_{22} = (T_2 ^{s} (t_s) - T^L (t_s))/P_2,
\end{aligned}
\end{equation}
The time constant $\theta_1=1/(C_T R_1 ^{e})$ is associated with thermal response of body 1 and computed using an effective resistance $R_1 ^{e}$, which is a weighted average of $R_{11}, R_{12}$. The thermal resistance terms can also be computed analytically for simple configurations as illustrated in Fig.~\ref{fig:2Body}. A general form of thermal conduction can be expressed as, $
R^{c} _{ij} = L_i/k_i A_i$. A set of trial runs conducted by varying the material properties and geometry reveal that a linear relationship exists between the empirically obtained and analytically computed thermal resistances $(R_{ij} = m_i R^{c}_{ij} + c_i)$ as presented in Fig.~\ref{fig:Rcalc}. \\
Now, for transient systems with power dissipation $P_1(t), P_2(t)$, a piecewise approximation is applied and the model equation is modified as,
\begin{equation}
    T^L (t)= \sum_j ^{n_T} \left[ \frac{P_{T,j}}{C_{T}} - \frac{P_{T,j-1}}{C_{T}} \right] (t-t^{0} _{j-1}),
\end{equation}
with $P_{T,j}=P_{1,j}+P_{2,j}$. This expression applies to any generalized system consisting of two or more bodies. The total power dissipation and thermal capacitance terms are modified as $P_{T,j}=\sum_i ^{N_B} P_{i,j}, C_T=\sum_i ^{N_B} C_i$, where $N_B$ represents the number of bodies in the system.
\begin{equation}
\begin{split}
    T^{D} _{1}(t) = \sum_j ^{n_T} [(P_{1,j} R_{11} + P_{2,j} R_{12}) \\ - (P_{1,j-1} R_{11} + P_{2,j-1} R_{12})] \\ (1-e^{-\theta_1 (t-t^{0} _{j-1})}), \\
    T^{D} _{2}(t) = \sum_j ^{n_T} [(P_{1,j} R_{21} + P_{2,j} R_{22}) \\ - (P_{1,j-1} R_{21} + P_{2,j-1} R_{22})] \\ (1-e^{-\theta_2 (t-t^{0} _{j-1})}),
\end{split}
\end{equation}
where, $j$ is the index of number of changes in the input power, $n_T$ is total number of changes (transients) in the input, and $t^{0} _{j}$ is the time instant at which the change in input powers occurs. The initial values are, $P_{1,0}=P_{2,0}=0, t^{0} _{0}=0$.
\subsection{Model for convection}
For convection heat transfer, with a single body, the problem inherently has two mediums: the solid and the surrounding fluid. The governing equation Eq.~\ref{eqmain1} must be modified to include the effects of flow velocity through continuity, Navier Stokes and energy equations, and the boundary condition must include convection boundary. For non-zero flow velocity, the temperature exhibits an exponential transient response, reaching a steady state more rapidly compared to purely conductive cases. This behavior is analogous to an association curve. The temperature evolution of the solid body in a fluid can be expressed as,
\begin{equation} \label{Eq:1BConv}
    T_1(t)=T^0 + T^L(t) + P_{1} R^{SF} _1 (1-e^{-\theta^{SF}_1 t}).
\end{equation}
The term $R^{SF} _1=1/hA$ denotes the thermal resistance between the solid body and the fluid, $h$ the convection heat transfer coefficient, $A$ the cross-sectional area, $\theta^{SF} _1=1/R^{SF}C^{SF}_T$ the time constant and $C^{SF} _{T}=C^S+C^F$ is the total thermal capacitance of the solid ($C^S$) and the fluid ($C^F$). The critical parameter here is the convection heat transfer coefficient that is often difficult to compute directly. However, this term can estimated empirically using the relation, 
\begin{equation} \label{Eq:h_1}
    h=\frac{Nu k}{D} = \frac{C_f Re^a Pr^b k}{D},
\end{equation}
where, $Nu$ is the Nusselt's number, $Re$ the Reynolds number, $Pr$ the Prandtl number, $k$ the thermal conductivity and $D$ the characteristic length. This can be further expanded as $Re=U_f D/\nu$ and $Pr=\nu/\alpha$ with $U_f$ the flow velocity, $\nu$ the viscosity, and $\alpha$ the thermal diffusivity. The coefficients $C_f,a,b$ are the model constants that depend on the geometry and flow conditions. Simplifying this we can express convection heat transfer as a function of geometry and flow properties as,
\begin{equation} \label{Eq:h_2}
    h = C_f U_f^a D^{a-1} \nu^{b-a} k^{1-b} \rho^b c_p^b.
\end{equation}
Typically, the model constants fall within the ranges of $a=0.5-1, b=0.3-0.4$. When $a>>b$ the convection heat transfer coefficient is primarily influenced by flow velocity $U_f$ and when $a \approx b$, $h$ is influenced by both fluid velocity and material properties. The convection coefficient is a field similar to the velocity and cannot be represented by a single lumped value for complex geometries. Accurate determination of $h$ requires multiple experimental runs and spatial averaging. However, to compute temperatures at specific points on the surface of the solid body, $h$ can be estimated empirically from the simulation data. In Ansys\textsuperscript{\textregistered} Icepak\textsuperscript{\texttrademark}, $h$ can be extracted from the option \textit{[Report $>$ Full Report {Variable: Convection Heat Transfer Coefficient}]}. A set of trial runs were performed at different flow velocities, where the $R^{SF}$ that best fits the characteristic temperature $T(t)$ for the given configuration, were estimated using curve-fitting. Convection coefficient and surface area were then extracted from Icepak\textsuperscript{\texttrademark} for the given configuration. Comparisons between the empirically estimated $h_{est}$ and the tool-reported $h_{calc}$ coefficients (Table~\ref{tab:h}) showed close agreement, confirming the validity of this modeling approach for the given configuration.
\begin{table}[h]
\begin{center}
\caption{Thermal resistance $R^{SF}$ estimated from curve-fitting $T(t)$, $A$ the surface area measured from the configuration. The tool reported convection heat transfer coefficient $h^{calc}$ is obtained from Ansys\textsuperscript{\textregistered} Icepak\textsuperscript{\texttrademark}.} 
\begin{tabular}{|c|c|c|c|c|}
\hline
\begin{tabular}[c]{@{}c@{}}$U^{f}$\\ $(m/s)$\end{tabular} & \begin{tabular}[c]{@{}c@{}}$R^{SF}=$\\ $1/h^{est}A$\\ $(K/W)$\end{tabular} & $A (m^2)$ & \begin{tabular}[c]{@{}c@{}}$h_{est}=$\\ $1/A R^{SF}$\\ $(W/m^2 K)$\end{tabular} & \begin{tabular}[c]{@{}c@{}}$h^{calc}$\\ $(W/m^2 K)$\end{tabular} \\ \hline
0.00098                                                      & 303.692                                                                     & 0.00015   & 21.95                                                                            & 21.83                                                            \\ \hline
5                                                            & 118.9                                                                       & 0.00015   & 56.07                                                                            & 56.19                                                            \\ \hline
50                                                           & 28.47                                                                       & 0.00015   & 234.09                                                                           & 235.7                                                            \\ \hline
\end{tabular} 
\label{tab:h}
\end{center}
\end{table}
For constant flow velocity or for flow velocities that vary by a small magnitude as in case of natural convection, a single value of $h$ can be used to estimate $R^{SF}$. However, when the flow velocity varies with time, $h(t)=f(U_f(t))$, it is observed from multiple simulations that for small range of velocities $U_f=U_0 \pm 10 m/s$ (with $U_0$ as a reference velocity), the relationship between thermal resistance, time constant and flow velocities can be expressed as, $	R_{ij} = m^R_i U_f ^{\gamma} + c^R_i, \theta_i = m^{\theta}_i U_f ^{\gamma} + c^{\theta}_i$.
Alternatively expressed, a linear relationship is observed between the parameters when the velocity is raised by a power $\gamma$. The model constants are obtained by running a parametric study for each heat source, for the specific range of velocities. For a two source system $P_1(t), P_2(t)$ with varying flow velocity $U_f(t)=[U_1(t),U_2(t),U_3(t)]$, where $U_1(t), U_2(t),U_3(t)$ are each a segment of velocity, we can develop a parametric study to generate the thermal resistance and time constants at a specific velocity.
\begin{enumerate}
	\item  $P_1(t)=1, P_2(t) = 0, U_f(t)=U_1$: 
	\begin{itemize} \item $R^{U_1} _{11}, R^{U_1} _{21}, \theta^{U_1} _1, \theta^{U_1} _2$
	\end{itemize}
	\item  $P_1(t)=1, P_2(t) = 0, U_f(t)=U_2$: 
	\begin{itemize} \item $R^{U_2} _{11}, R^{U_2} _{21}, \theta^{U_2} _1, \theta^{U_2} _2$
	\end{itemize}
	\item  $P_1(t)=1, P_2(t) = 0, U_f(t)=U_3$: 
	\begin{itemize} \item $R^{U_3} _{11}, R^{U_3} _{21}, \theta^{U_3} _1, \theta^{U_3} _2$
	\end{itemize}
	
	\item  $P_1(t)=0, P_2(t) = 1, U_f(t)=U_1$: 
	\begin{itemize} \item $R^{U_1} _{12}, R^{U_1} _{22}, \theta^{U_1} _1, \theta^{U_1} _2$
	\end{itemize}
	\item  $P_1(t)=0, P_2(t) = 1, U_f(t)=U_2$: 
	\begin{itemize} \item $R^{U_2} _{12}, R^{U_2} _{22}, \theta^{U_2} _1, \theta^{U_2} _2$
	\end{itemize}
	\item  $P_1(t)=0, P_2(t) = 1, U_f(t)=U_3$: 
	\begin{itemize} \item $R^{U_3} _{12}, R^{U_3} _{22}, \theta^{U_3} _1, \theta^{U_3} _2$
	\end{itemize}
\end{enumerate}
Now, for all the velocities the relationship between the variables are obtained as,
\begin{equation}
	\begin{aligned}
	R_{ij} (U) = m^{R} _i U_f ^{\gamma_k} + c^{R} _i, \\
	\theta_i (U) = m^{\theta}_i U_f ^{\gamma_l} + c^{\theta} _i, \\
	\end{aligned}
\end{equation}
Here, the model constants $m^R, c^R, m^{\theta}, c^{\theta}, \gamma_k, \gamma_l$ are empirically derived by minimizing the error between thermal resistance and time constants from simulation and those estimated from modeling. Finally, the temperature of each component is determined by applying superposition considering the temperature rise due to applied power dissipation, variable velocity and effects of additional heat sources. This is expressed as,
\begin{equation}
	\begin{aligned}
	T_1 (t) = T^0 + T^L(t) + \sum_j ^{n_T} \big[ P_1(t) R_{11} (U) [1-e^{-\theta_1 (U) (t-t^{0} _{j-1})}] \\ + P_2(t) R_{12} (U) [1-e^{-\theta_2 (U) (t-t^{0} _{j-1})}] \big], \\
	T_2 (t) = T^0 + T^L(t) + \sum_j ^{n_T} \big[ P_1(t) R_{21} (U) [1-e^{-\theta_1 (U) (t-t^{0} _{j-1})}] \\ + P_2(t) R_{22} (U) [1-e^{-\theta_2 (U) (t-t^{0} _{j-1})}] \big],
\end{aligned}
\end{equation}
where, $j$ is the index of number of changes in the input power, $n_T$ is total number of changes (transients) in the input, and $t^{0} _{j}$ is the time instant at which the change in input powers occurs.
\subsection{Estimation of model constants for the case study of inverter module with mixed modes of heat transfer}
For general multibody systems involving both conduction and convection heat transfer modes as in Fig.~\ref{fig:T_comp}, the model constants are estimated by the following process, 
\begin{enumerate}
	\item \label{param1} A parametric CFD simulation is performed by applying 1W of input power to one MOSFET junction while keeping the inputs of the other MOSFETs as 0 W. The simulation is run for a short duration $(t_s \approx 10 s, \Delta t \approx 0.2 s$), since the power dissipation for highly transient power electronics systems typically vary at a rate of $100-1000 \mu s$. The temperatures of all six MOSFETs (sources), the PCBA and the heatsink are recorded. This process is repeated for each individual heat source. Thus, the number of trials equals the number of sources present in the system. 
	\item For each trial run, the characteristic thermal resistances and time constants for each component are determined by minimizing the sum of squared errors $\bigl[SSE =\sum_i \bigl((T^{s} _i (t)-T^L(t)) - T_i(t) \bigr)^2 \bigr]$ between the temperature difference measured from the simulation and the linear temperature curve $T^{s}_i-T^L(t)$ and the temperature predicted by the model $T_i(t)$. A generalized reduced gradient optimization method is used to estimate the parameters that minimize the $SSE$. 
	\item In Trial 1, the effect of MOSFET 1 on all other MOSFETs, the PCBA, and the heatsink is obtained. For example, $R_{11}$ represents the thermal resistance measured at MOSFET 1 when power is applied at MOSFET 1, while $R_{31}$ represents the thermal resistance measured at MOSFET 3 under the same input condition. In trial 2, the influence of MOSFET 2 on all other components is computed, and so on. 
	\item To compute the temperature of MOSFET 1 considering all heat sources, the superposition of contributions from each source is used as,
	\begin{equation}
		\begin{aligned}
			T_1 = T^0 + T^L(t) + P_1(t) R_{11} [1-e^{-\theta_{1} t}] + P_2(t) R_{12} [1-e^{-\theta_{2} t}] + \\ 
			P_3(t) R_{13} [1-e^{-\theta_{3} t}] +  P_4(t) R_{14} [1-e^{-\theta_{4} t}] + \\
			P_5(t) R_{15} [1-e^{-\theta_{5} t}] +  P_6(t) R_{16} [1-e^{-\theta_{6} t}].
			\end{aligned}
	\end{equation}
	A similar approach is used for all the other components. 
	
	\item For forced convection systems with variable flow velocity, the parametric simulations described in Step~\ref{param1} must be extended to account for variations in flow velocity. The modified parametric simulations are defined as,
	\begin{itemize}
		\item $P_1 (t)=1, P_2(t), \ldots, P_6(t)=0, U_f=1$ 
		\item $P_1 (t)=1, P_2(t), \ldots, P_6(t)=0, U_f=5$
		\item $P_1 (t)=1, P_2(t), \ldots, P_6(t)=0, U_f=10$
		\item $P_1 (t)=1, P_2(t), \ldots, P_6(t)=0, U_f=20$
	\end{itemize}
	The same set of velocity conditions is applied for each heat source. Therefore, for 6 sources and 4 velocity magnitudes, a total of 24 simulations are required.
\end{enumerate}
\section{ROM for LTI systems}
LTI-ROM is primarily applicable to systems that are both linear and time-invariant. The thermal system presented in this work can be modeled as an LTI system in the case of forced convection with constant flow velocity. The system can be considered linear if the total input power dissipation can be represented as the sum of contributions from the power dissipated by each individual component. Similarly, the flow can be considered linear if the total flow velocity induced (in the presence of multiple flow sources) is equal to the sum of the flows induced by each source individually. This system can be considered time invariant if the fluid properties like density, viscosity, flow rate, and operating conditions remain constant over time. However, if non-linear effects like variable flow velocity, transition to turbulence due to non-linear interactions or variable operating conditions such as density driven flow as in natural convection are introduced, the LTI assumption breaks down. \\
For multiple input multiple output (MIMO), LTI system with 6 power sources and 8 temperature monitor points (including PCBA, HS) as presented in Sec.~\ref{casestudy}, the state-space representation is,
\begin{equation}
	\begin{aligned}
		\dot{x} = A x(t) + B u(t), \\
		y(t) = C x(t)  + D.
	\end{aligned}
\end{equation}
In this formulation, 
\begin{itemize}
\item $x=x_{ij} (t) \in \mathbb{R}^{48}$ (for 6 sources, 8 monitor points), 
\item $A=-k_{ij} \in \mathbb{R}^{48 \times 48}$ (time constant matrix), 
\item $B=1/k_{ij} R_{ij} \in \mathbb{R}^{48 \times 6}$ (thermal capacitance matrix), 
\item $u(t) = P(t) \in \mathbb{R}^{6}$ (input power dissipation), 
\item $y(t) = T(t) \in \mathbb{R}^{8}$ (output temperature), 
\item $C \in \mathbb{R}^{8 \times 48}$ (summation of effect of all the sources), 
\item $D = T_0$ (ambient temperature). \\
\end{itemize}
The corresponding temperature transfer function is,
\begin{equation}
	T_i (t) = T_0 +\sum_{i=1} ^{6} (\psi_{ij}  \ast  P_i)(t),
\end{equation}
where, $\psi_{ij}$ is the impulse function. This is computed by taking a derivative of the temperature function and expressed as,
\begin{equation}
	\psi_{ij} = k_{ij} R_{ij} e^{-k_{ij} t}.
\end{equation}
For discrete systems, the temperature for new inputs are computed as,
\begin{equation}
	T_i [n] = T_0 + \sum_{i=1} ^{6} (\psi_{ij}  \ast  P_i)[n].
\end{equation}
For LTI systems, the thermal resistance and time constants are computed by characterizing the temperature response of the system after reaching steady state. The LTI assumption requires linearity around a baseline and this steady state is typically where the system is assumed to behave linearly. While there are similarities in assumptions and modeling strategies between LTI-ROM and LPLSP, the key difference lies in the duration over which the output temperature is characterized. LTI-ROM requires modeling the system’s response until steady-state temperature is achieved, whereas LPLSP models transient temperature behavior over shorter durations. The applicability of this method to forced convection with constant velocity is straightforward since the model coefficients are constants. However, in the case of forced convection with variable velocity, an interpolated or parameter varying ROM will be required. This extension is beyond the scope of this work. Therefore, comparisons between LPLSP and LTI-ROM in this study are limited to cases of natural convection and forced convection with constant flow velocity.
\section{Results}
\begin{figure}
	\centering
	\begin{subfigure}[b]{0.45\textwidth}
		\centering
	\includegraphics[width=\linewidth]{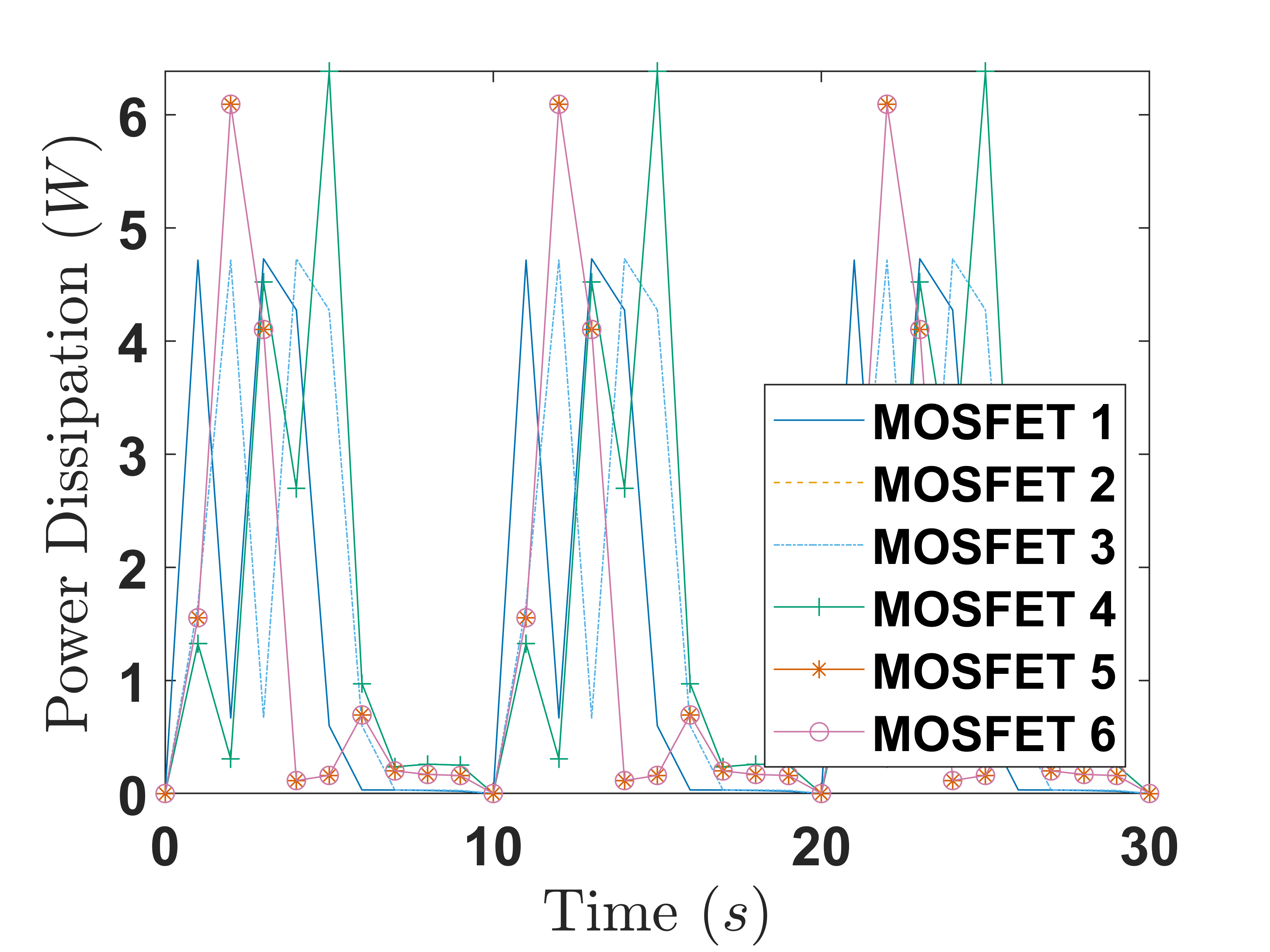}
	\subcaption{Input power dissipation}  
	\label{fig:power}
\end{subfigure} \hfill
\begin{subfigure}[b]{0.45\textwidth}
	\centering
	\includegraphics[width=\linewidth]{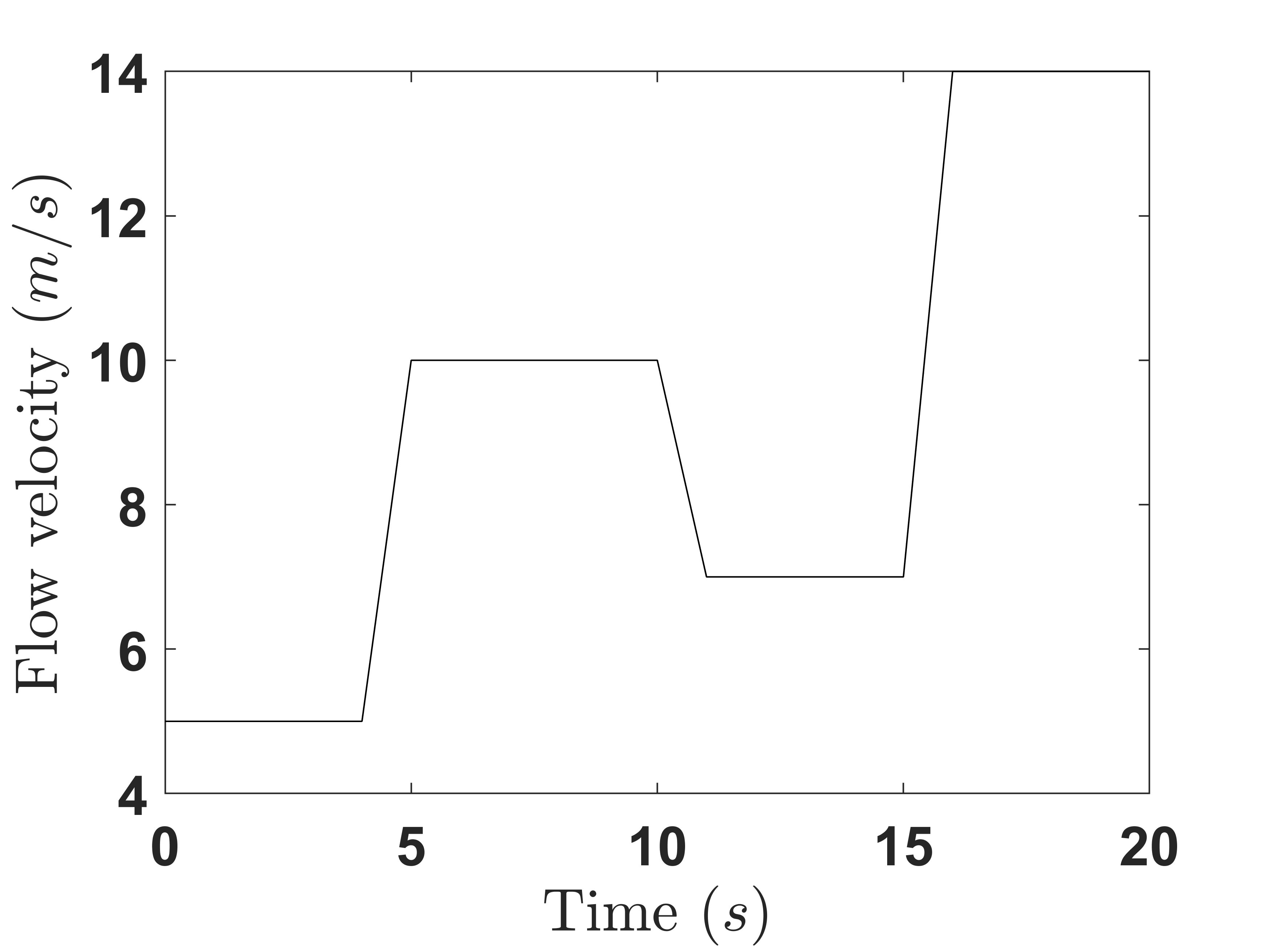}
	\subcaption{Velocity profile} \label{fig:velocity}
\end{subfigure}
	\caption{Same input power dissipation is applied for the cases of natural convection, forced convection with constant velocity and forced convection with variable flow velocity. Variable flow velocity profile $u_x (m/s)$ is applied only to the case of forced convection with variable flow velocity. This simulation is run only for a duration of $20 s$.}
	\label{fig:input}
\end{figure}
\begin{figure}
	\centering
	\begin{subfigure}[b]{0.48\textwidth}
		\centering
		\includegraphics[width=\linewidth]{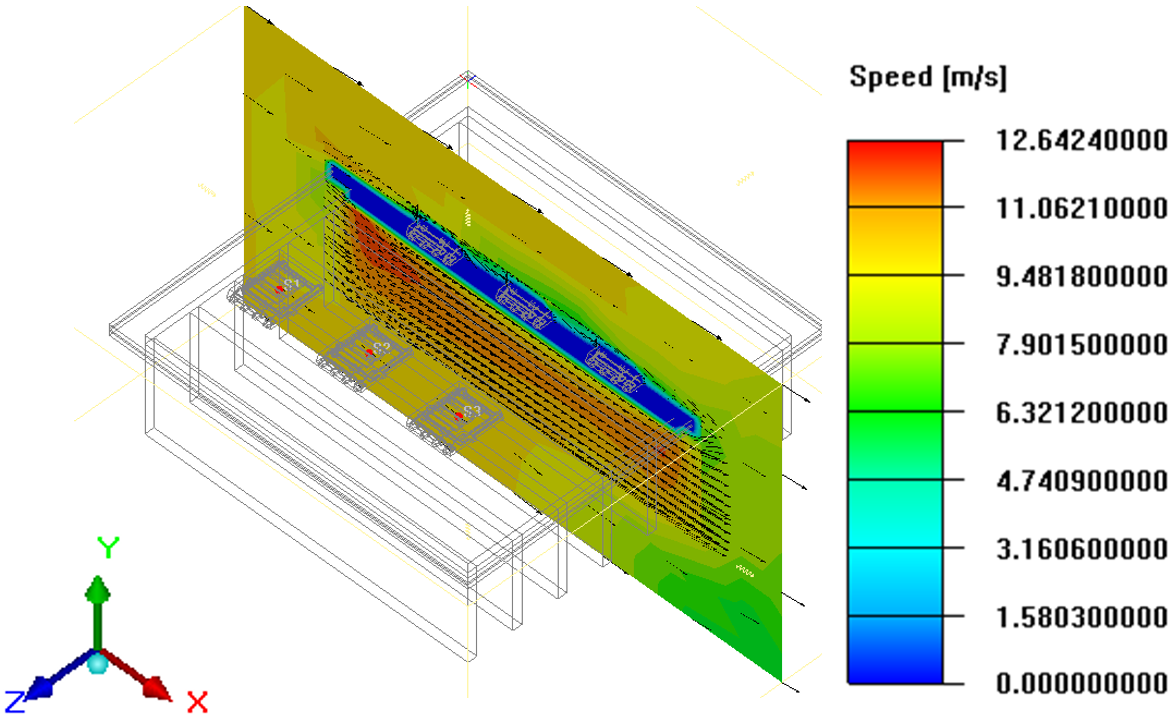}
		\subcaption{Contour plot of speed $U=\sqrt{u_x^2 + u_y^2 + u_z^2}$ (Isometric view)} \label{fig:speed} 
	\end{subfigure} \hfill
	\begin{subfigure}[b]{0.48\textwidth}
		\centering
		\includegraphics[width=\linewidth]{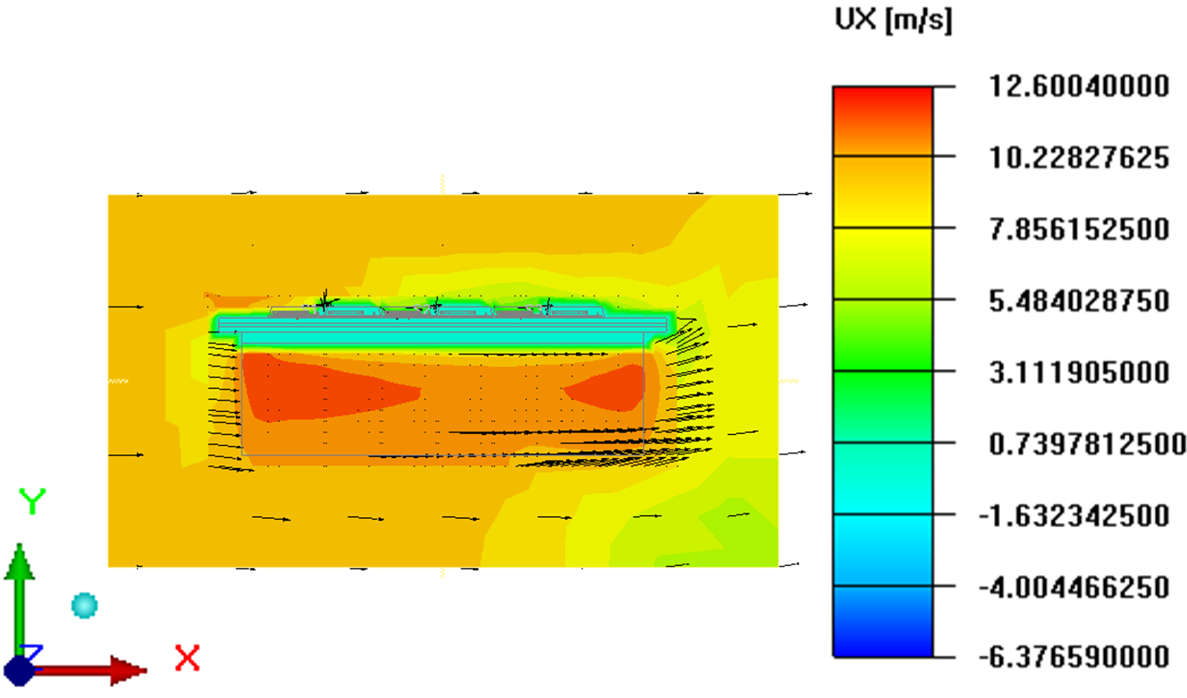}
		\subcaption{Contour plot of velocity $u_x$ viewed from $z$ axis} \label{fig:velocity_contour_z}
	\end{subfigure} \\
	\begin{subfigure}[b]{0.48\textwidth}
		\centering
		\includegraphics[width=\linewidth]{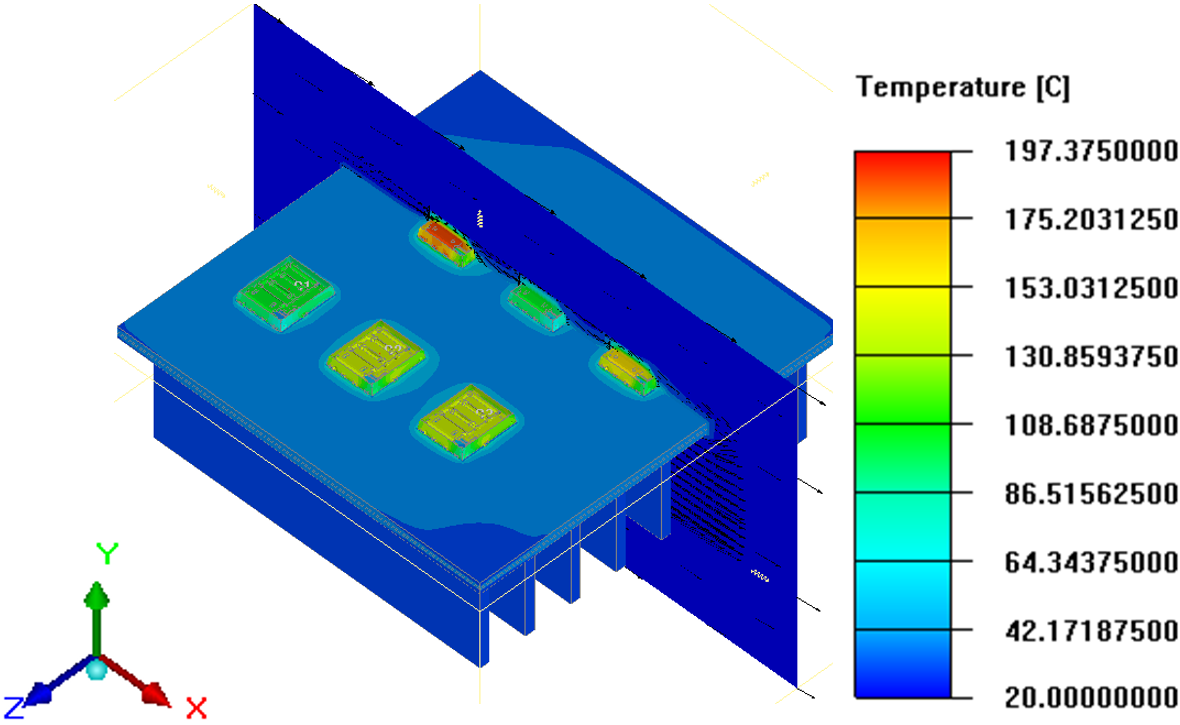}
		\subcaption{Contour plot of Temperature (Isometric view)} \label{fig:temperature_contour_iso}
	\end{subfigure} \hfill
		\begin{subfigure}[b]{0.48\textwidth}
		\centering
		\includegraphics[width=\linewidth]{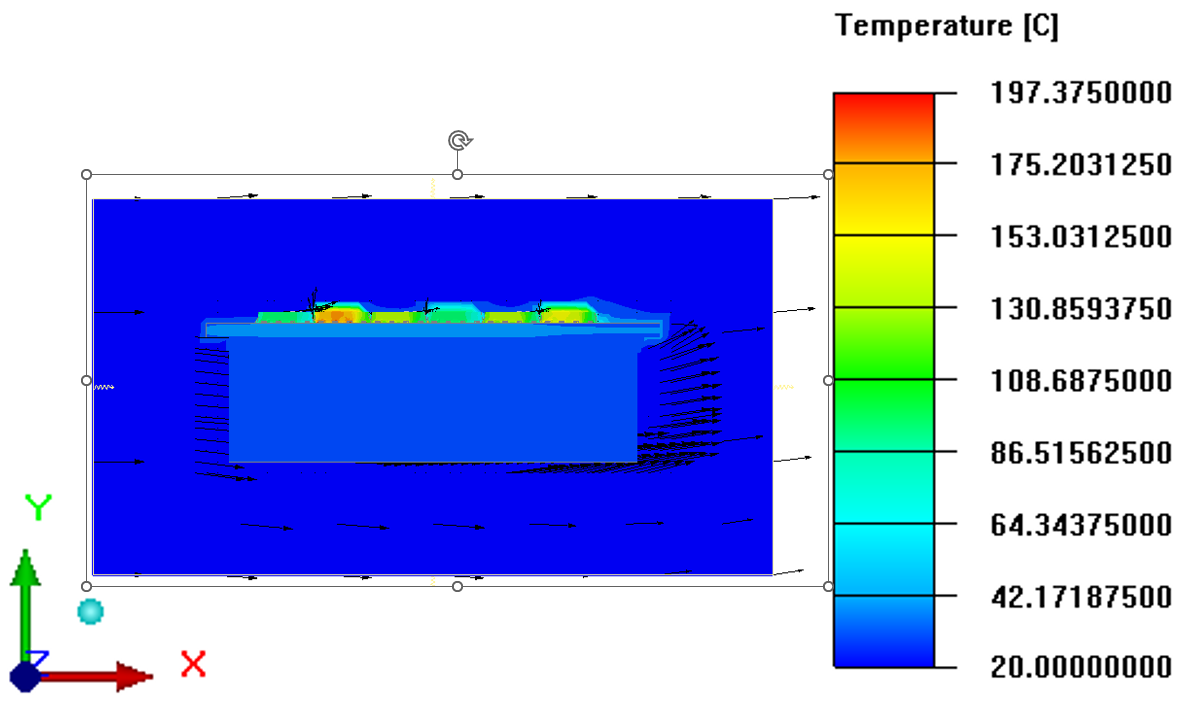}
		\subcaption{Contour plot of Temperature viewed from $z$ axis} \label{fig:temperature_contour_z}
		\end{subfigure}
	\caption{Contour plots of velocity magnitudes and surface temperatures from CFD simulation in Ansys\textsuperscript{\textregistered} Icepak\textsuperscript{\texttrademark} at simulation time of $t=20 s$, for the case of forced convection with variable flow velocity.}
	\label{fig:input}
\end{figure}

\begin{figure}[htbp]
	\centering
	\begin{subfigure}[b]{0.45\textwidth}
		\centering
		\includegraphics[width=\linewidth]{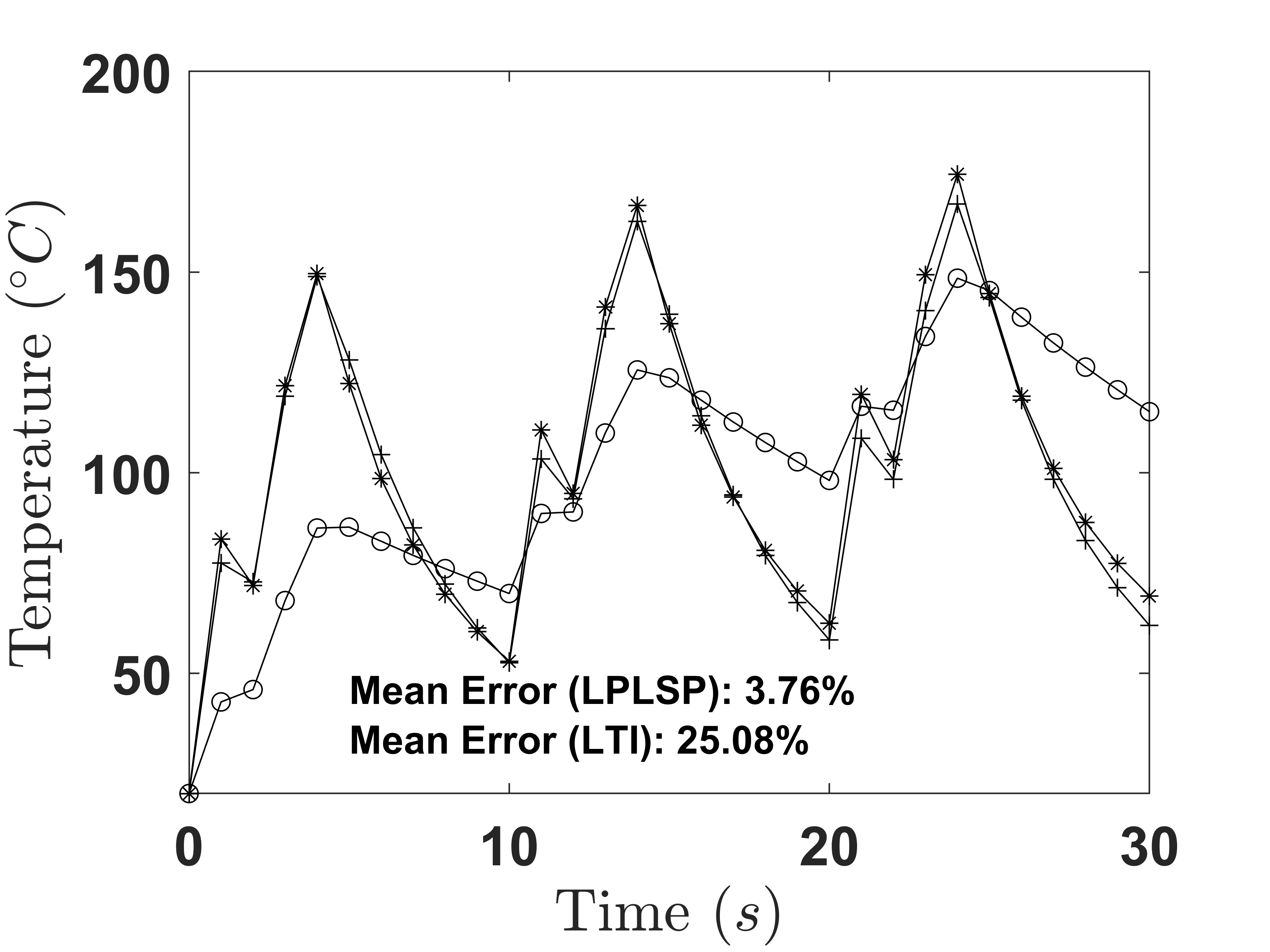}
		\subcaption{MOSFET 1}
	\end{subfigure}
	\begin{subfigure}[b]{0.45\textwidth}
		\centering
		\includegraphics[width=\linewidth]{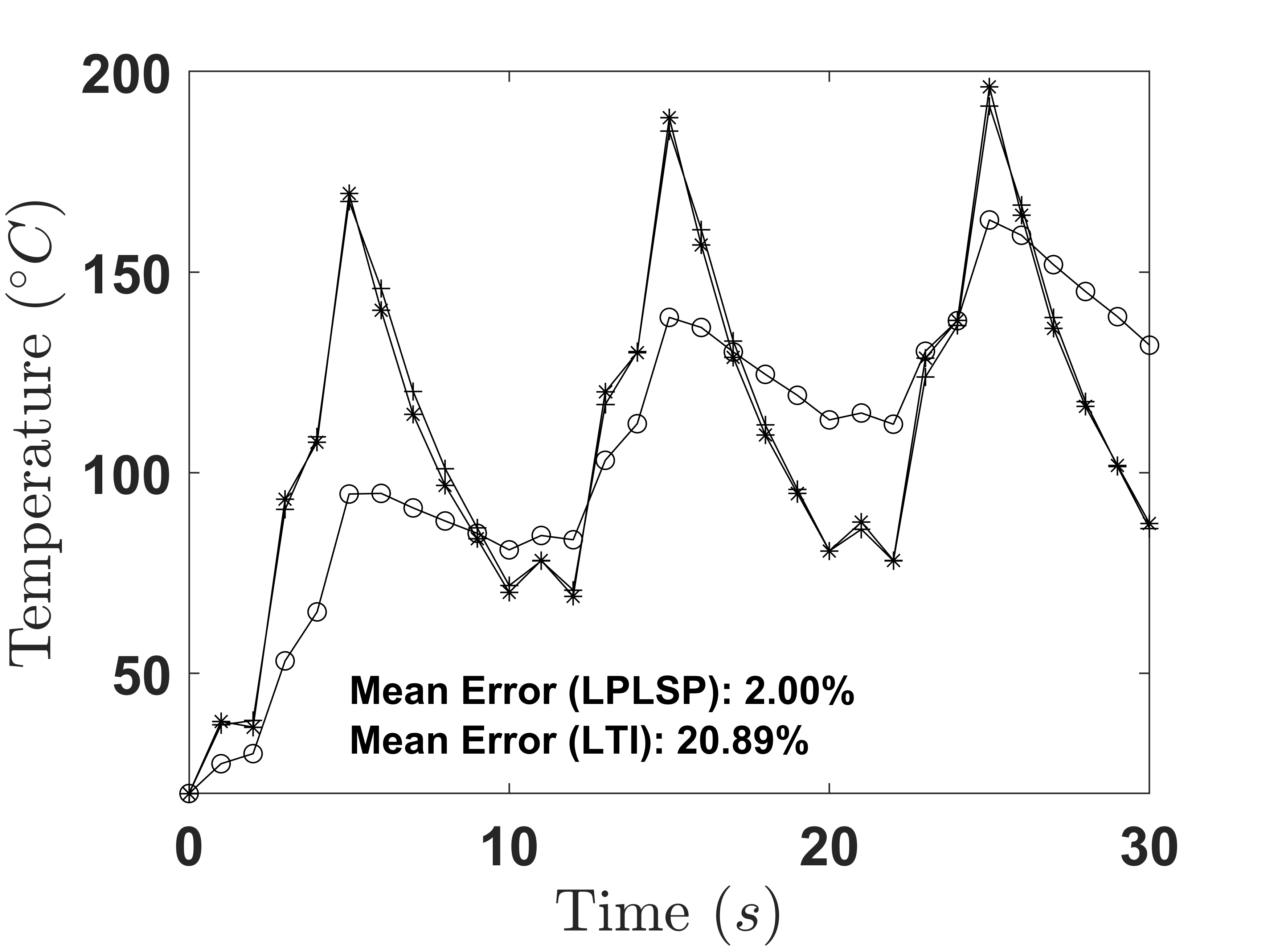}
		\subcaption{MOSFET 2}
	\end{subfigure}
	\begin{subfigure}[b]{0.45\textwidth}
		\centering
		\includegraphics[width=\linewidth]{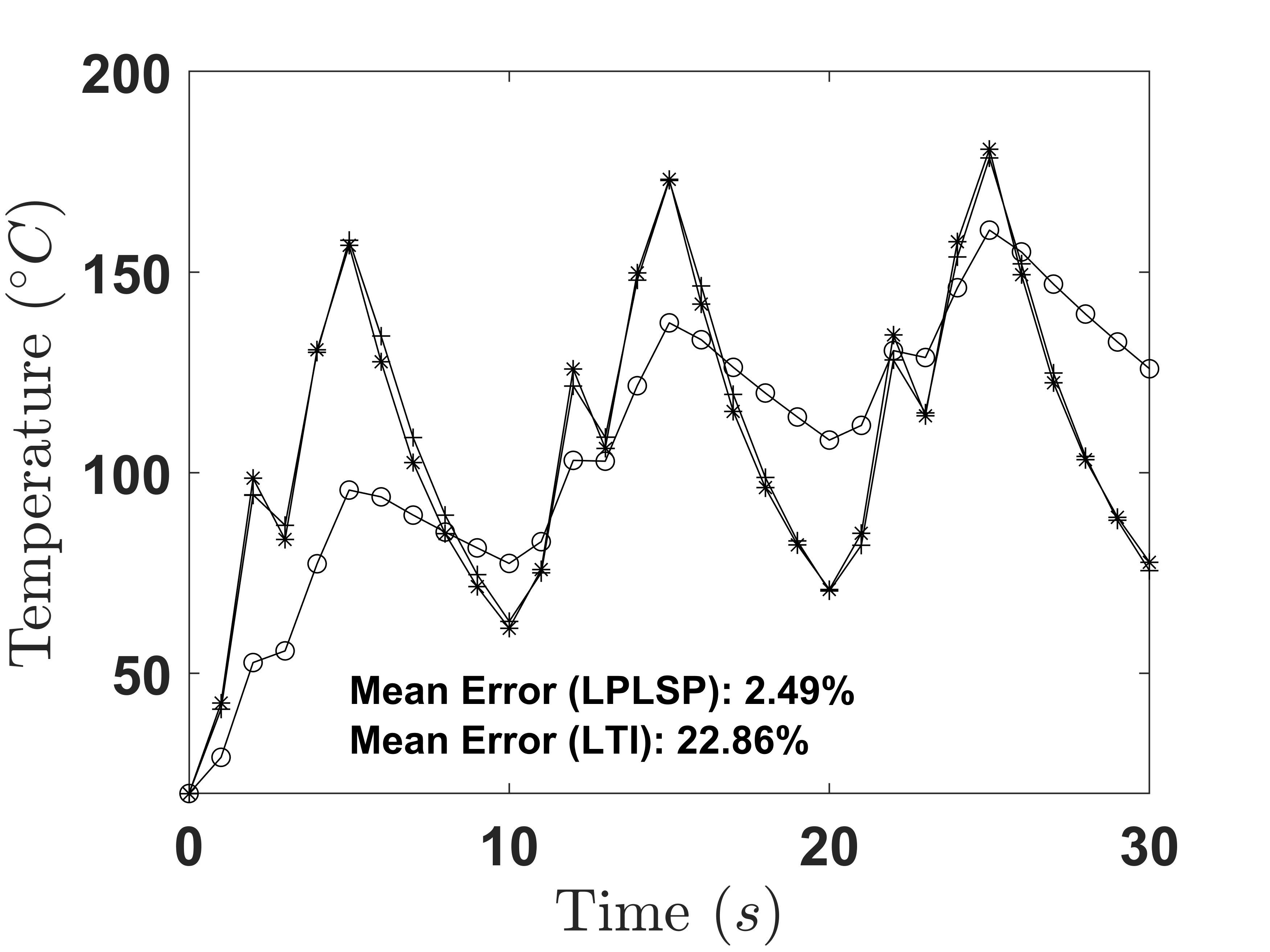}
		\subcaption{MOSFET 3}
	\end{subfigure}
	\begin{subfigure}[b]{0.45\textwidth}
		\centering
		\includegraphics[width=\linewidth]{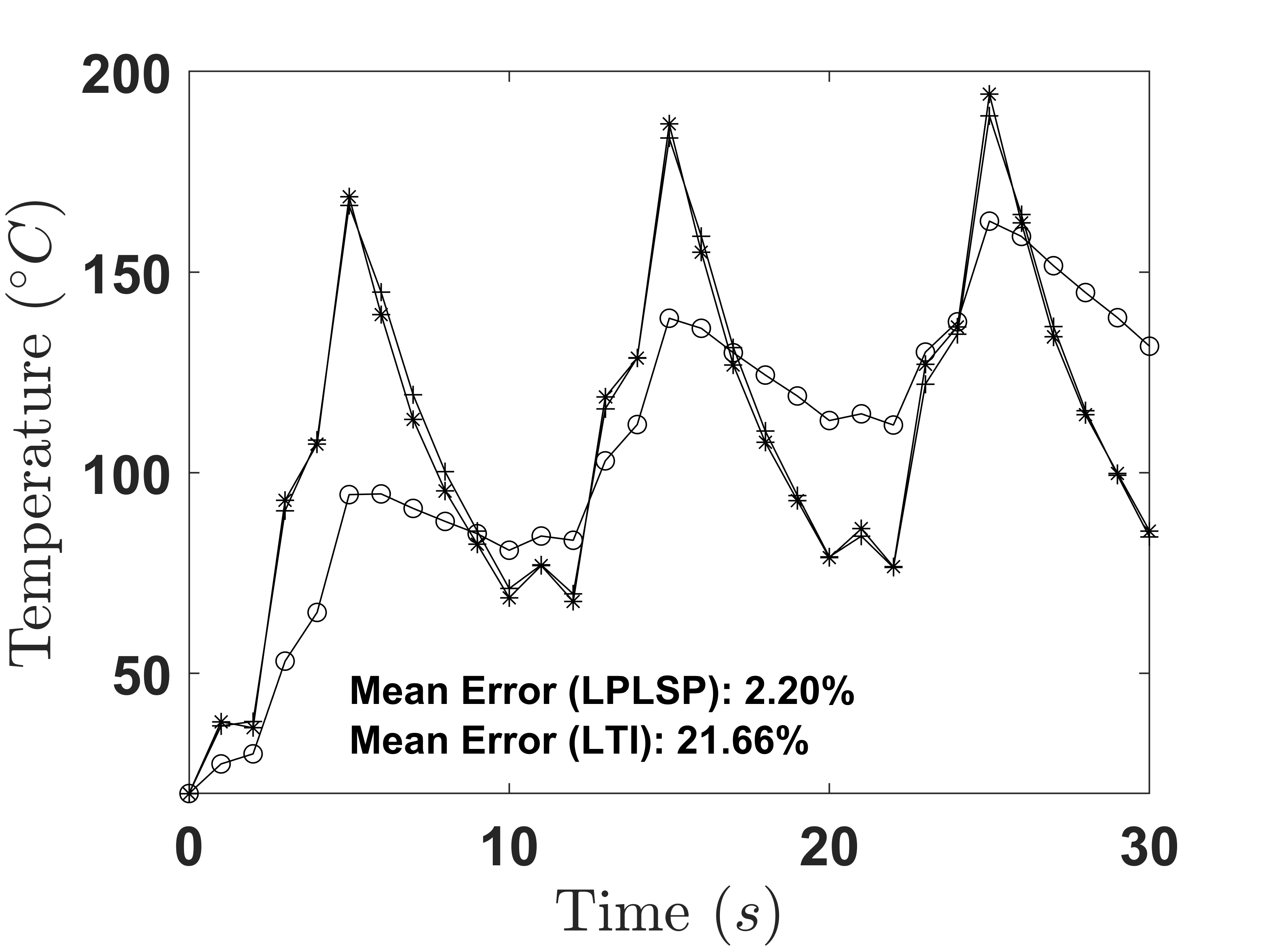}
		\subcaption{MOSFET 4}
	\end{subfigure}
	\begin{subfigure}[b]{0.45\textwidth}
	\centering
	\includegraphics[width=\linewidth]{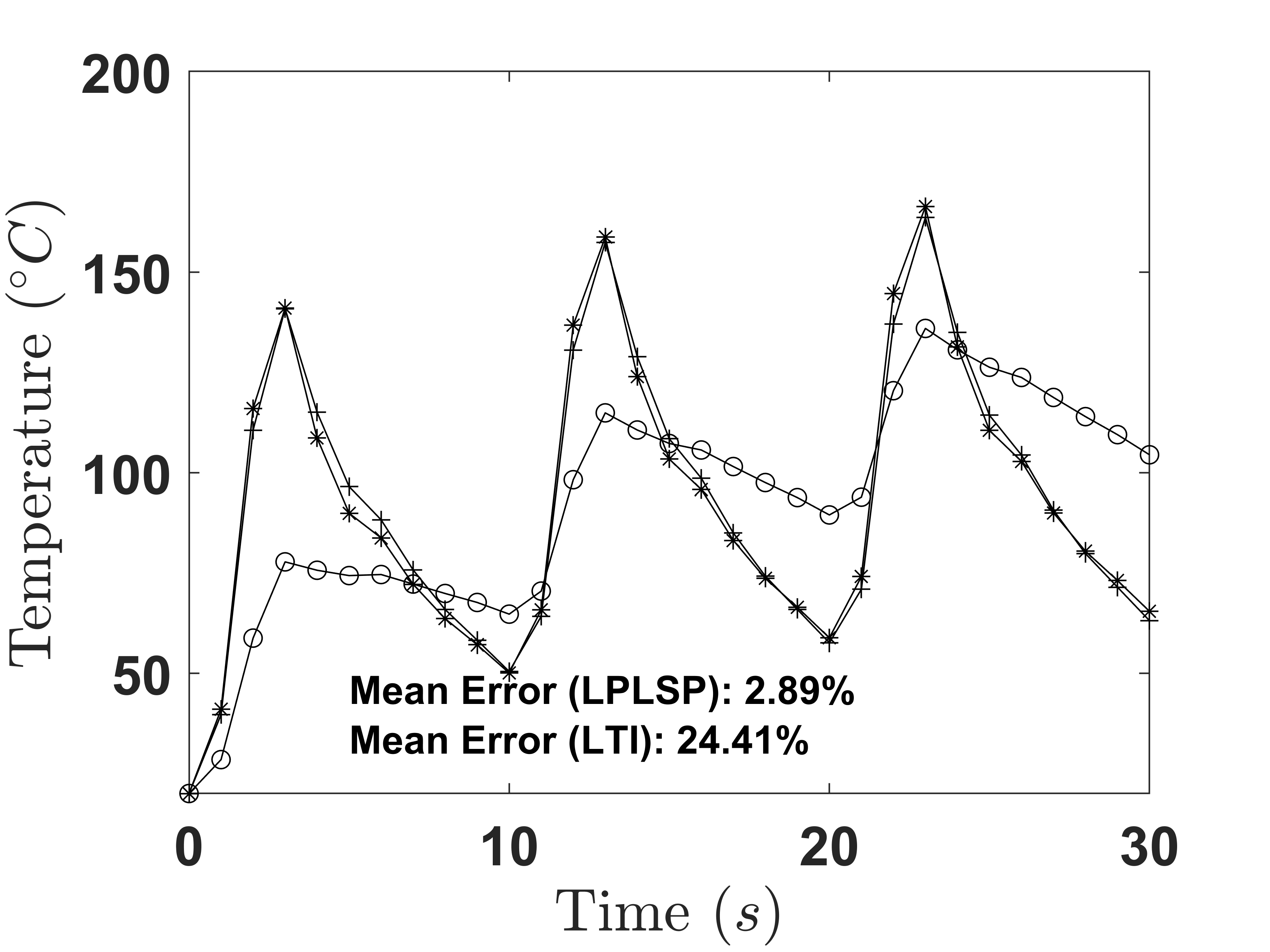}
	\subcaption{MOSFET 5}
	\end{subfigure}
	\begin{subfigure}[b]{0.45\textwidth}
		\centering
		\includegraphics[width=\linewidth]{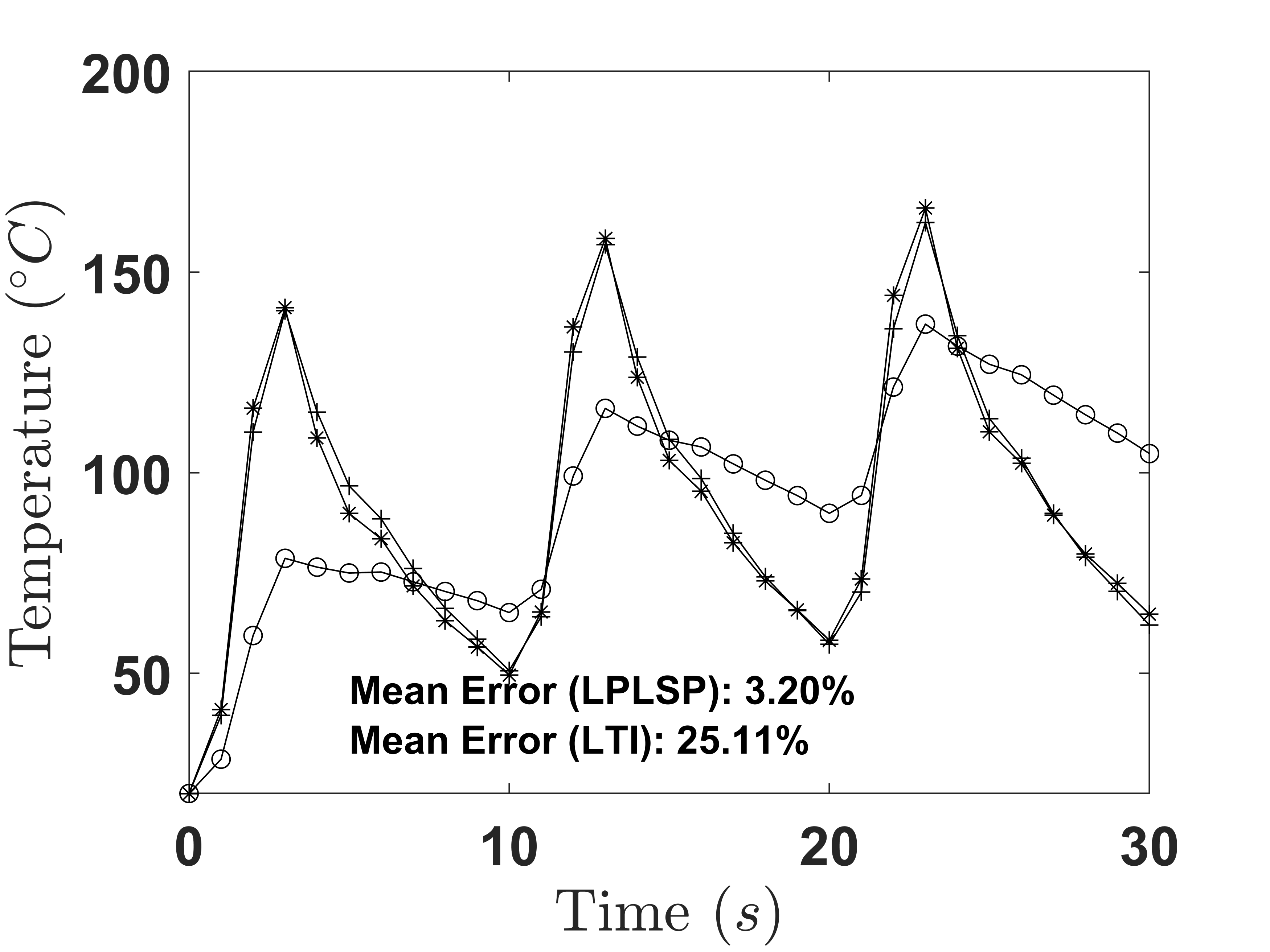}
		\subcaption{MOSFET 6}
	\end{subfigure}	
	\caption{Comparison of temperatures from simulation $(-*)$, LPLSP model $(-+)$ and LTI-ROM $(-o)$ for the case of inverter module in a natural convection environment. Mean percentage error $\overline{Error}\% = |\frac{T^{s}_i-T_i}{\max(T^{s}_i)}|$ between the model and simulation data is also presented.}
	\label{fig:NC}
\end{figure}
\begin{figure}[htbp]\ContinuedFloat
	\centering
	\begin{subfigure}[b]{0.45\textwidth}
		\centering
		\includegraphics[width=\linewidth]{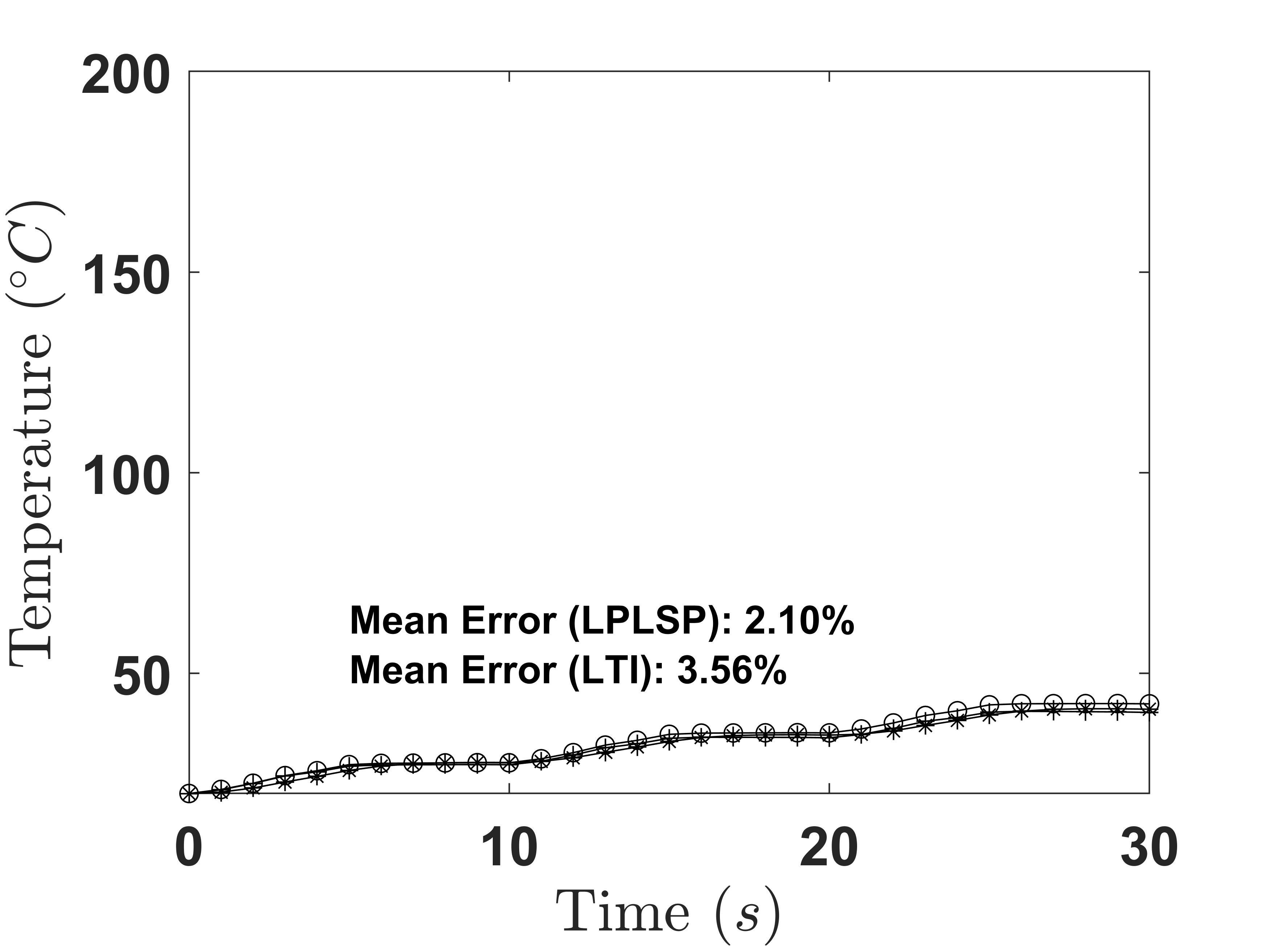}
		\subcaption{PCBA}
	\end{subfigure}
	\begin{subfigure}[b]{0.45\textwidth}
		\centering
		\includegraphics[width=\linewidth]{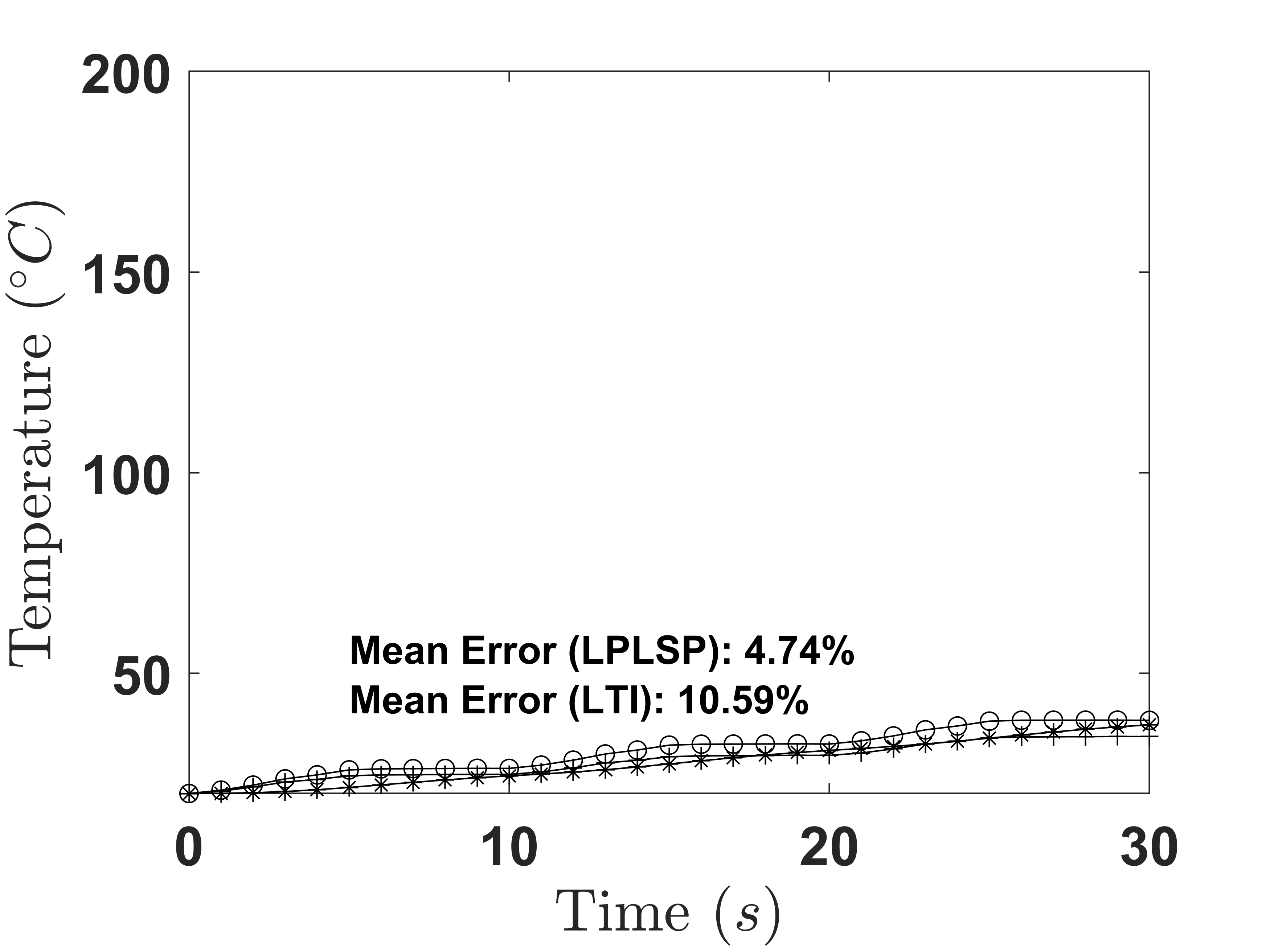}
		\subcaption{Heat Sink}
	\end{subfigure}
	\caption*{(Figure \ref{fig:NC} continued)}
\end{figure}
\begin{figure}
	\centering
	\begin{subfigure}[b]{0.45\textwidth}
	\centering
	\includegraphics[width=\linewidth]{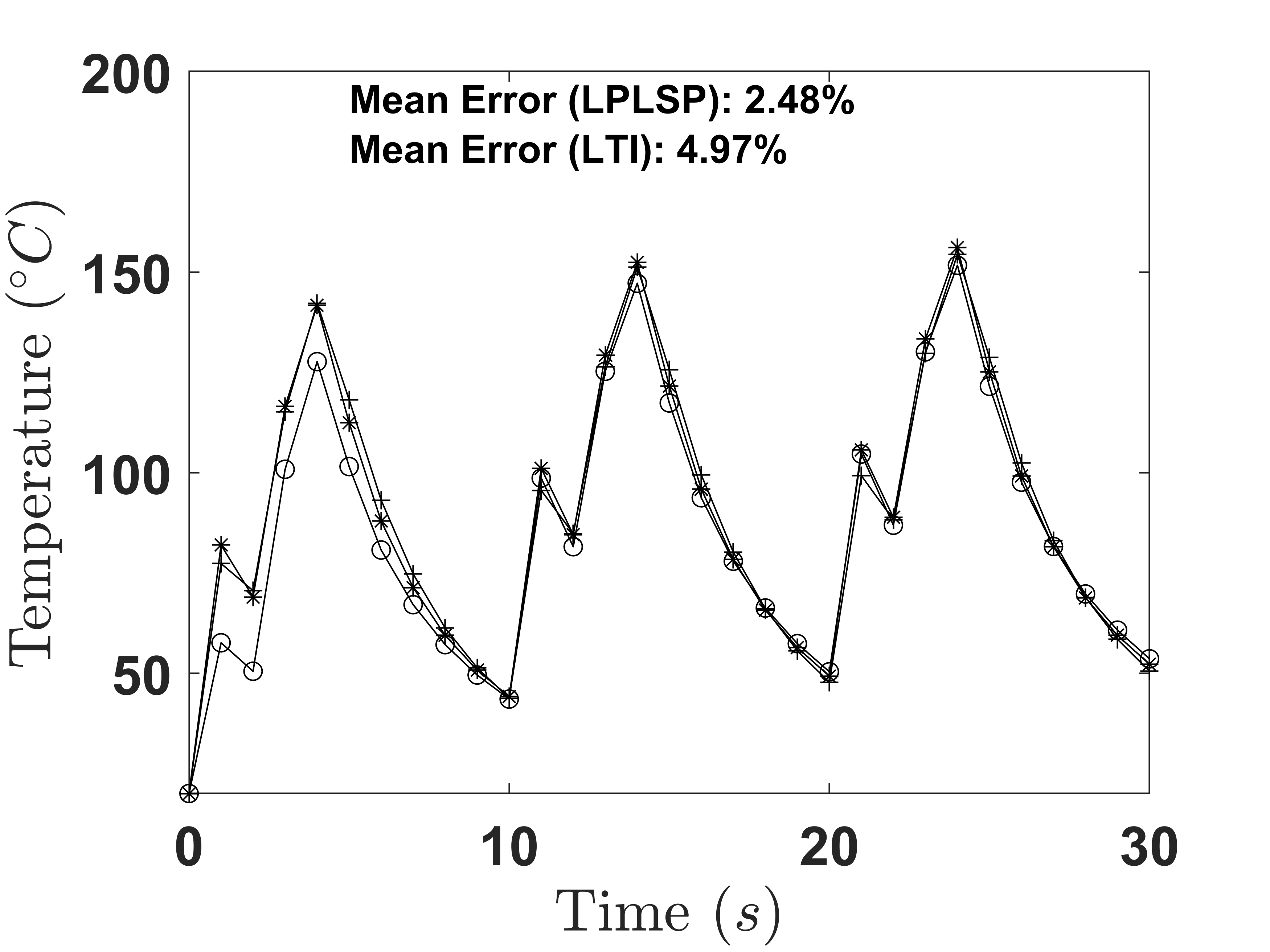}
	\subcaption{MOSFET 1}  
\end{subfigure} 
\begin{subfigure}[b]{0.45\textwidth}
	\centering
	\includegraphics[width=\linewidth]{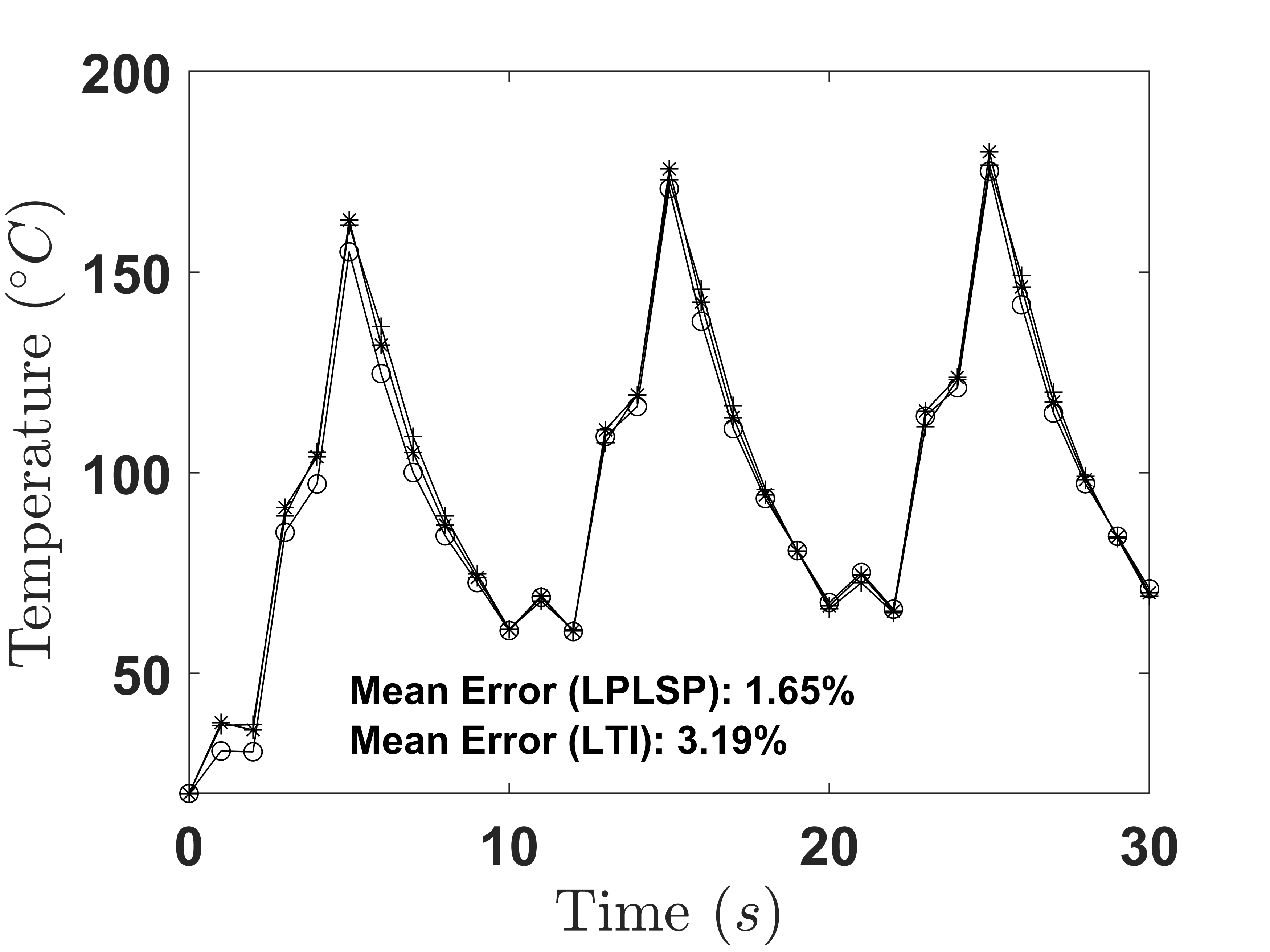}
	\subcaption{MOSFET 2} 
\end{subfigure} 
\begin{subfigure}[b]{0.45\textwidth}
	\centering
	\includegraphics[width=\linewidth]{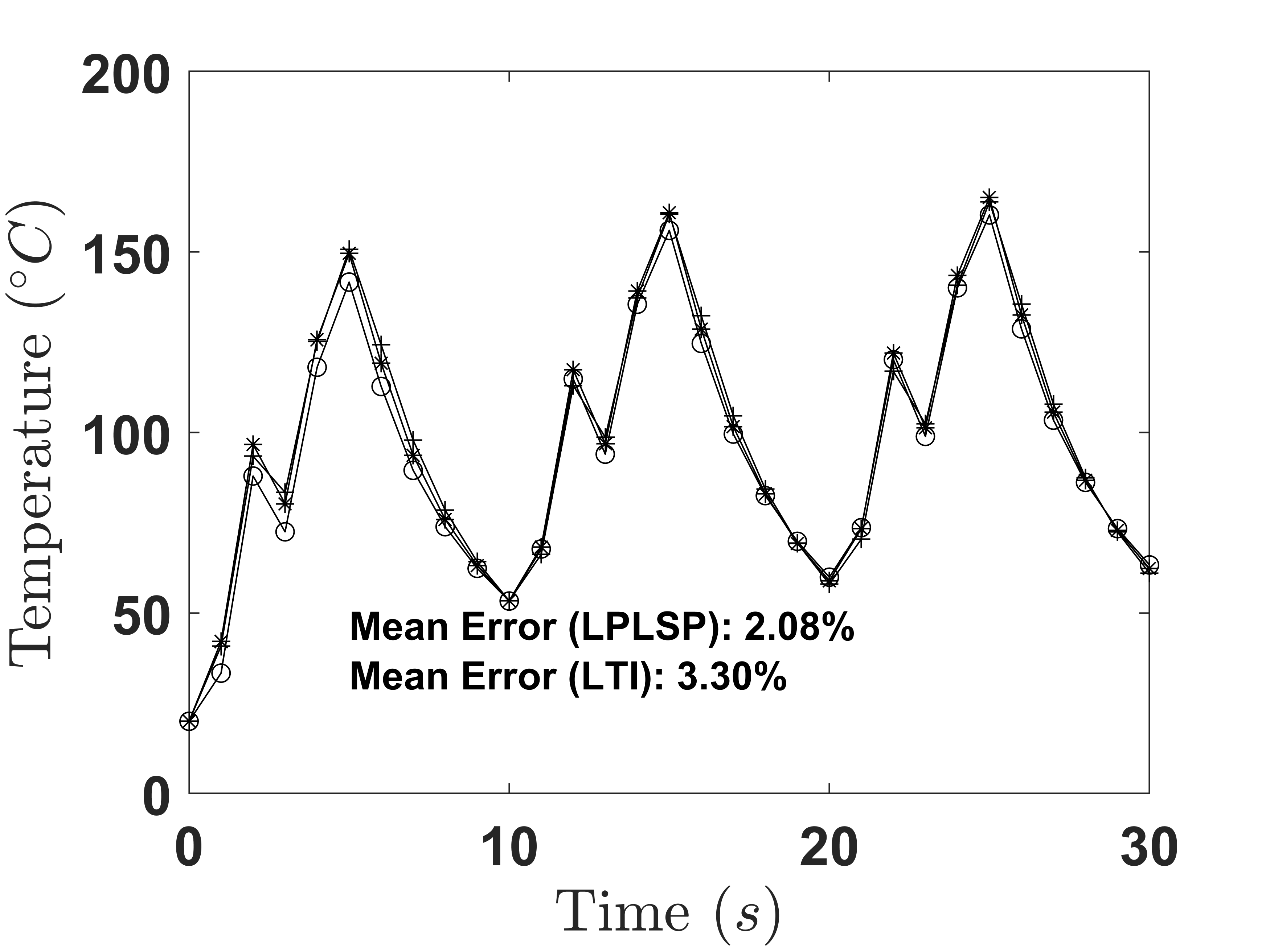}
	\subcaption{MOSFET 3}  
\end{subfigure} 
\begin{subfigure}[b]{0.45\textwidth}
	\centering
	\includegraphics[width=\linewidth]{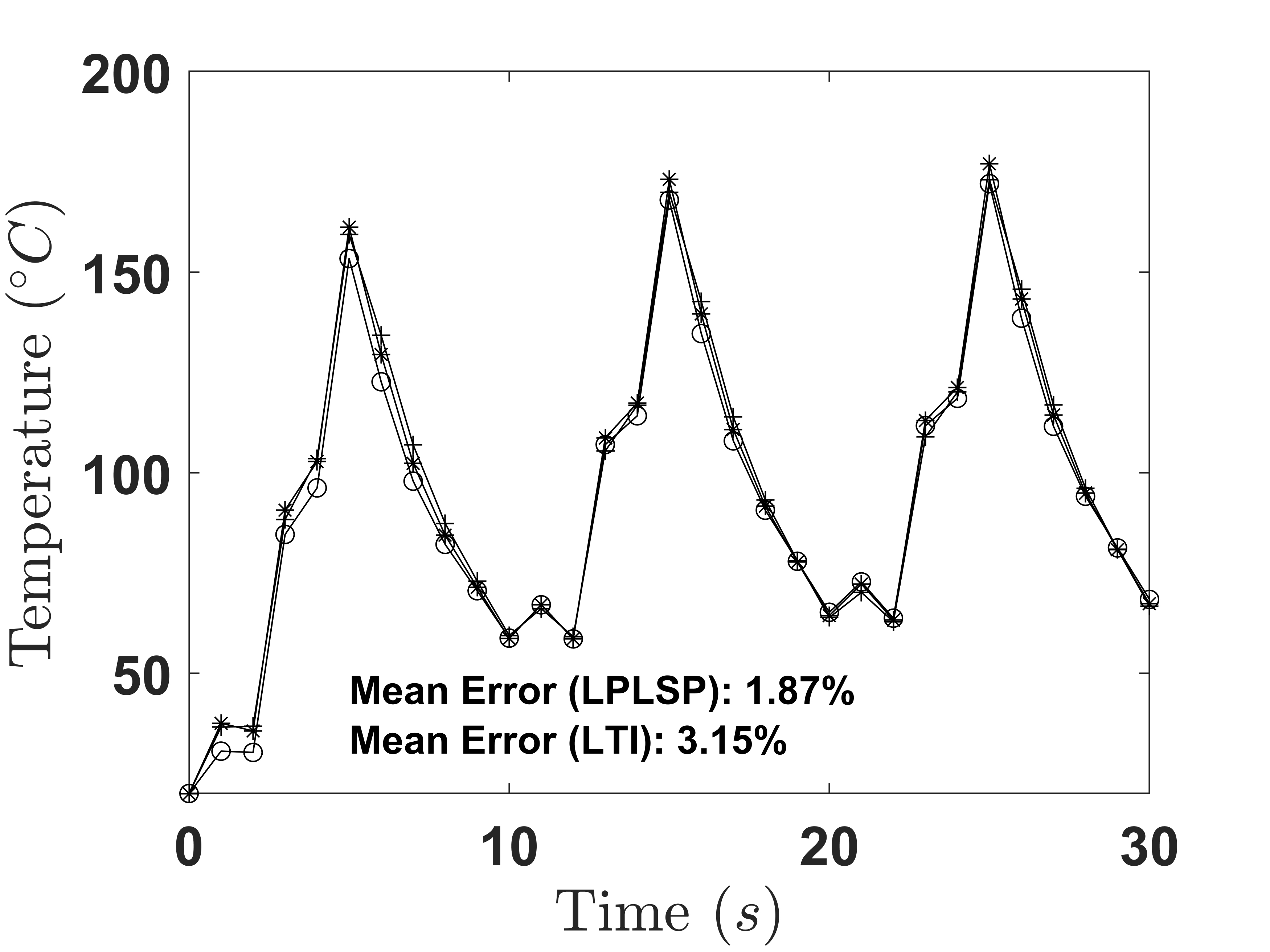}
	\subcaption{MOSFET 4} 
\end{subfigure} 
\begin{subfigure}[b]{0.45\textwidth}
	\centering
	\includegraphics[width=\linewidth]{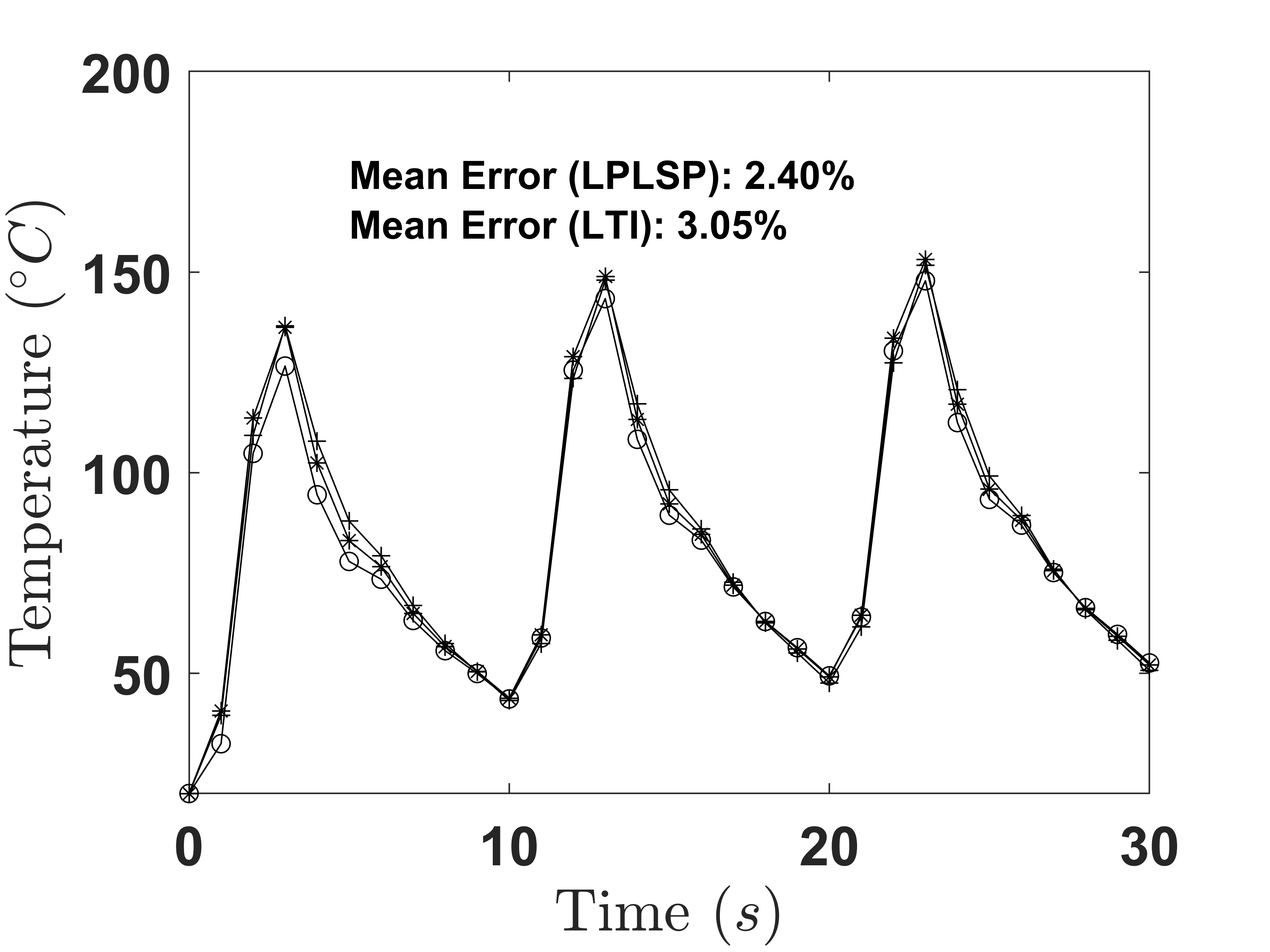}
	\subcaption{MOSFET 5}  
\end{subfigure}
\begin{subfigure}[b]{0.45\textwidth}
	\centering
	\includegraphics[width=\linewidth]{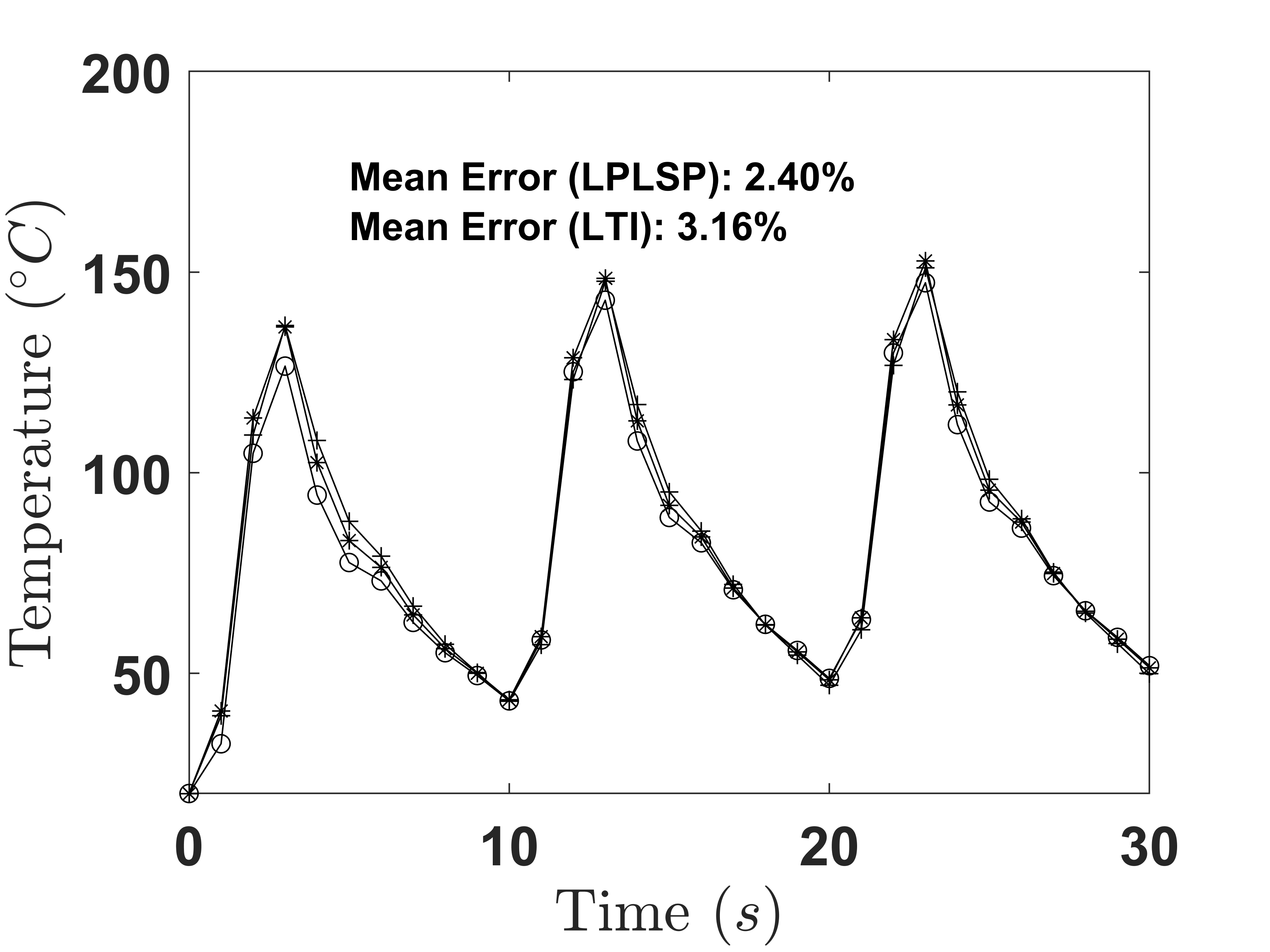}
	\subcaption{MOSFET 6} 
\end{subfigure} 
	\caption{Comparison of temperatures from simulation $(-*)$, LPLSP model $(-+)$ and LTI-ROM $(-o)$ for the case of inverter module in a forced convection environment with a constant flow velocity of $u_x=5 m/s$. Mean percentage error $\overline{Error}\% = |\frac{T^{s}_i-T_i}{\max(T^{s}_i)}|$ between the model and simulation data is also presented.}
\label{fig:FC_const}
\end{figure}
\begin{figure}[htbp]\ContinuedFloat
	\centering
\begin{subfigure}[b]{0.45\textwidth}
	\centering
	\includegraphics[width=\linewidth]{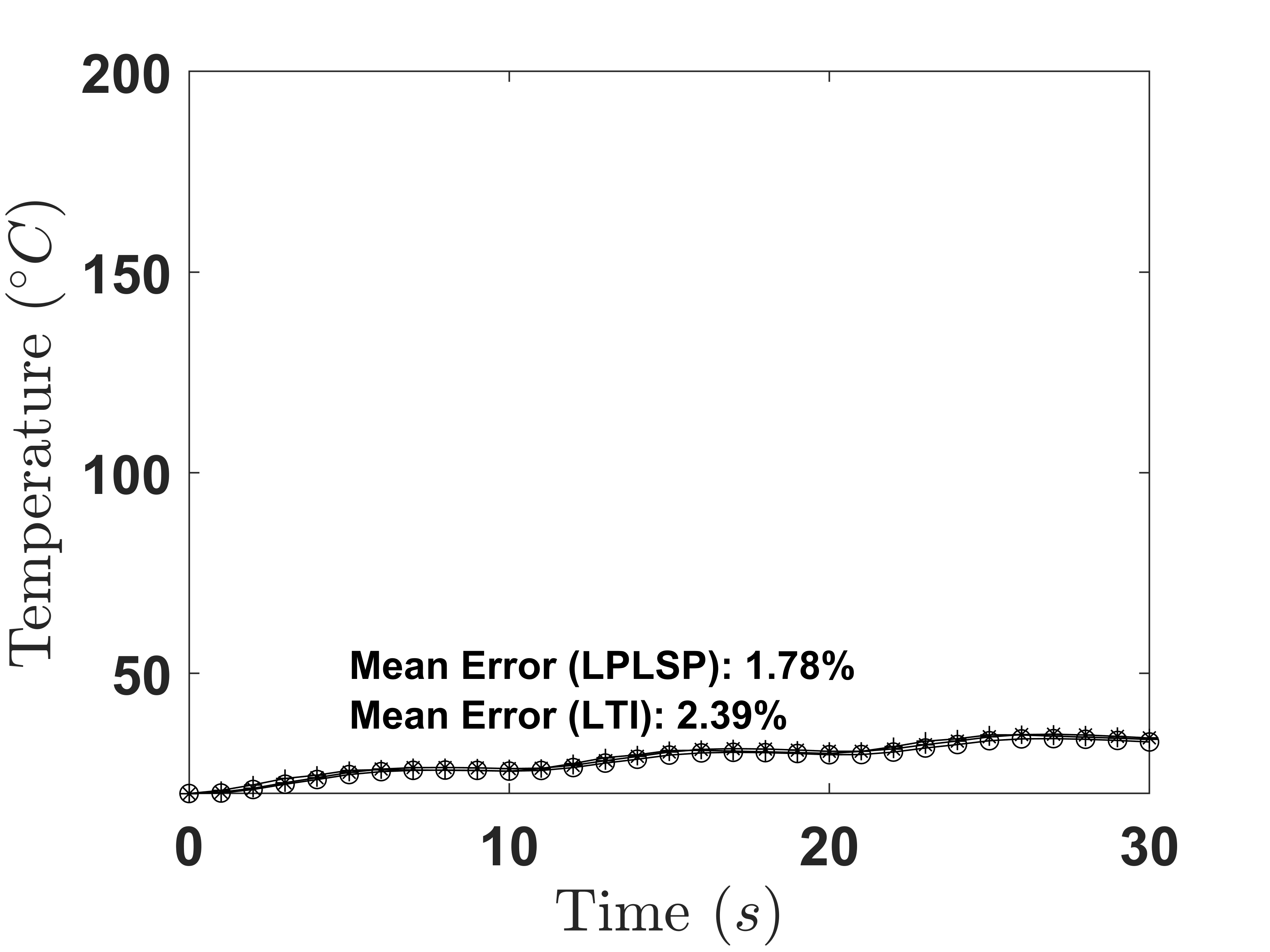}
	\subcaption{PCBA}  
\end{subfigure} 
\begin{subfigure}[b]{0.45\textwidth}
	\centering
	\includegraphics[width=\linewidth]{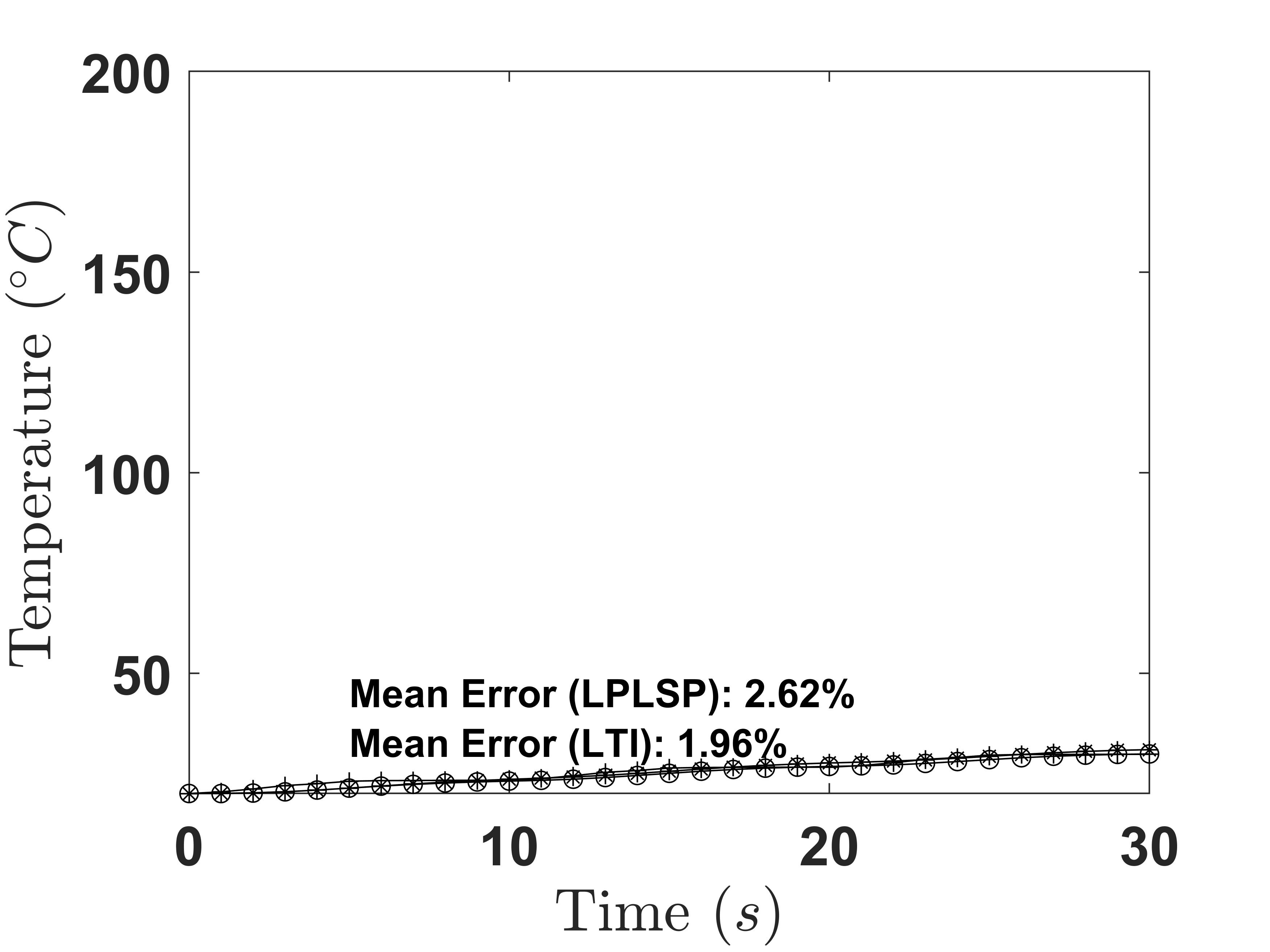}
	\subcaption{Heat Sink} 
\end{subfigure}
	\caption*{(Figure \ref{fig:FC_const} continued)}
\end{figure}
\begin{figure}
	\centering
	\begin{subfigure}[b]{0.45\textwidth}
		\centering
		\includegraphics[width=\linewidth]{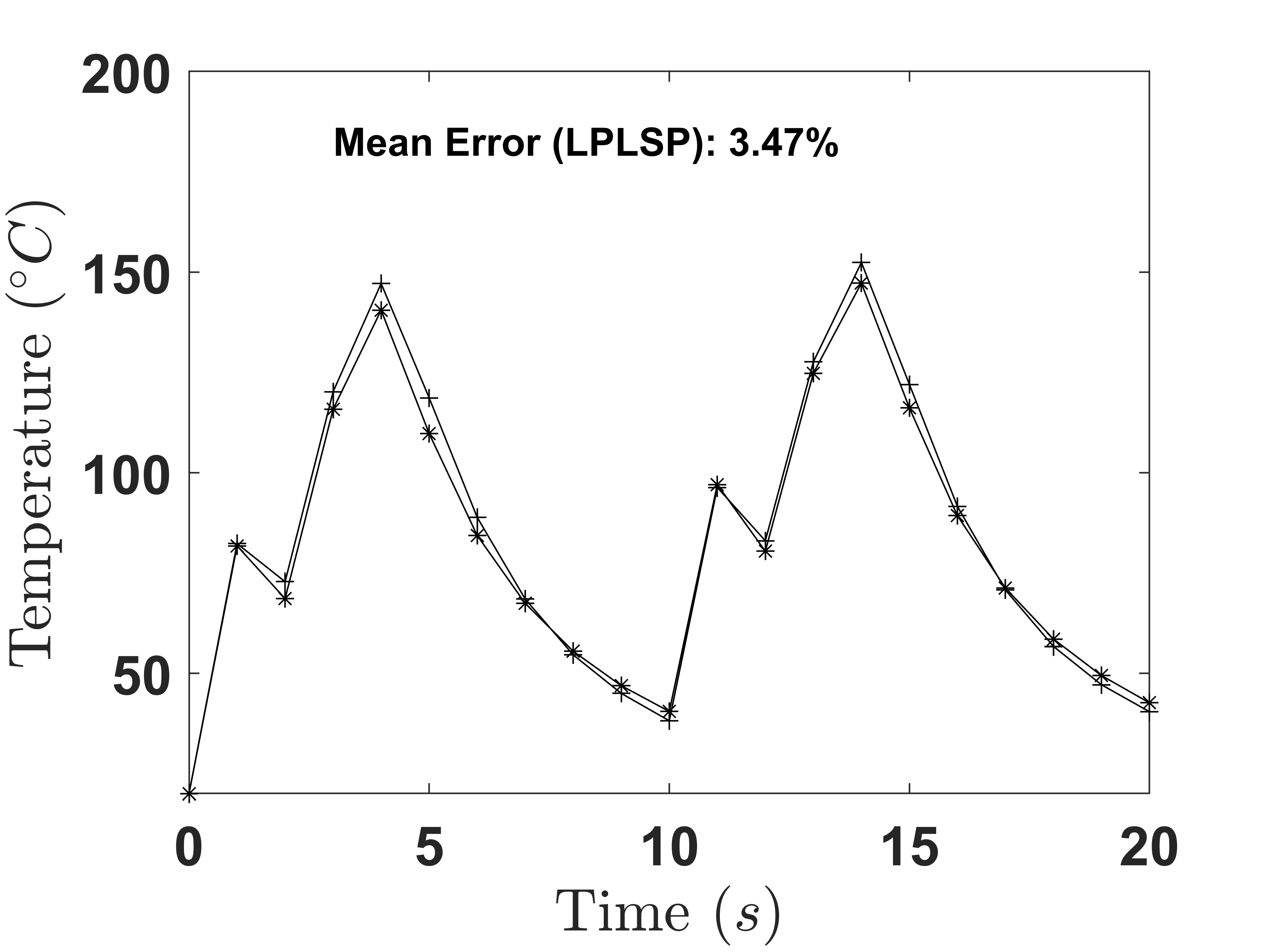}
		\subcaption{MOSFET 1}  
	\end{subfigure} 
	\begin{subfigure}[b]{0.45\textwidth}
		\centering
		\includegraphics[width=\linewidth]{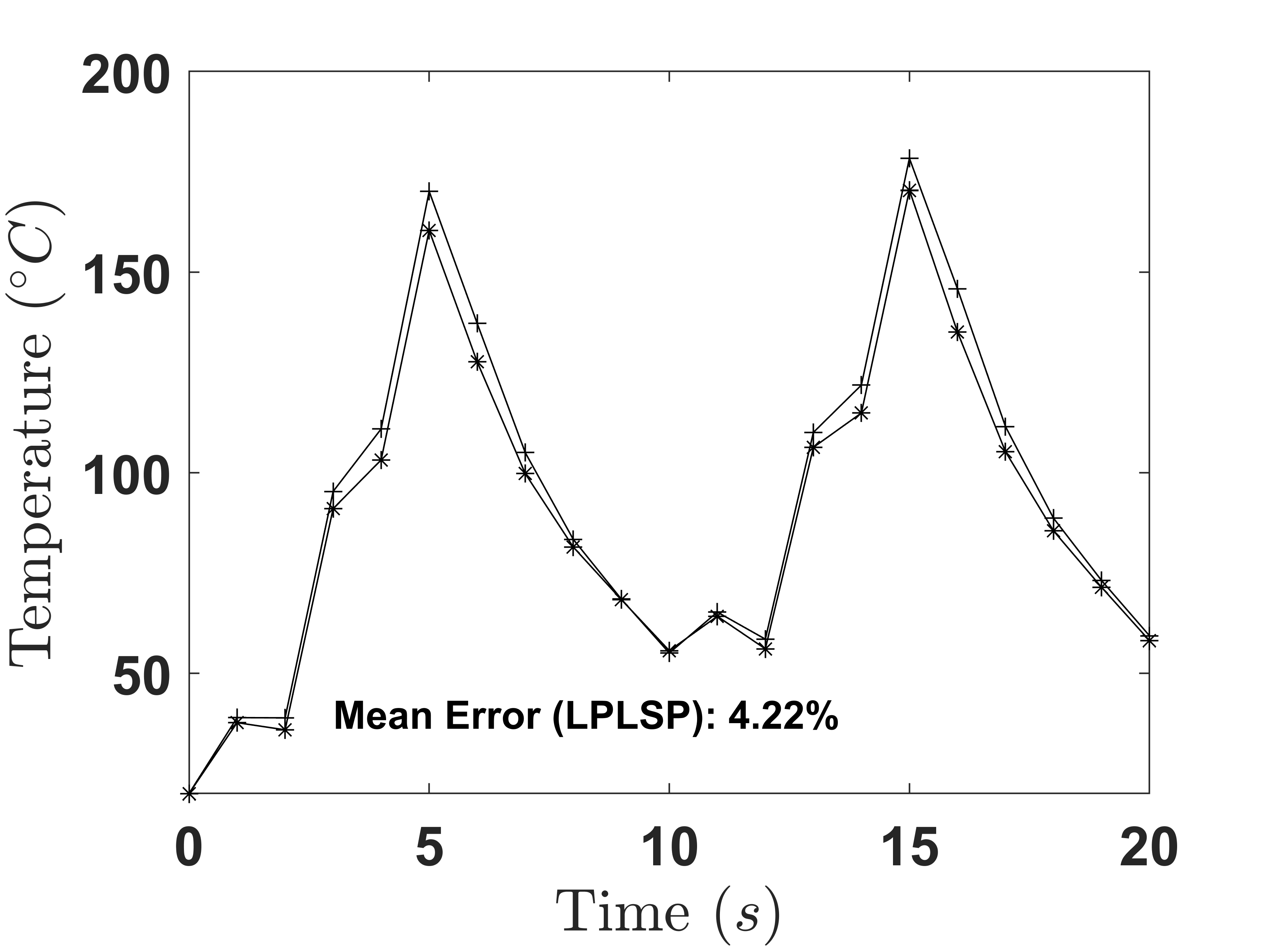}
		\subcaption{MOSFET 2} 
	\end{subfigure} 
	\begin{subfigure}[b]{0.45\textwidth}
		\centering
		\includegraphics[width=\linewidth]{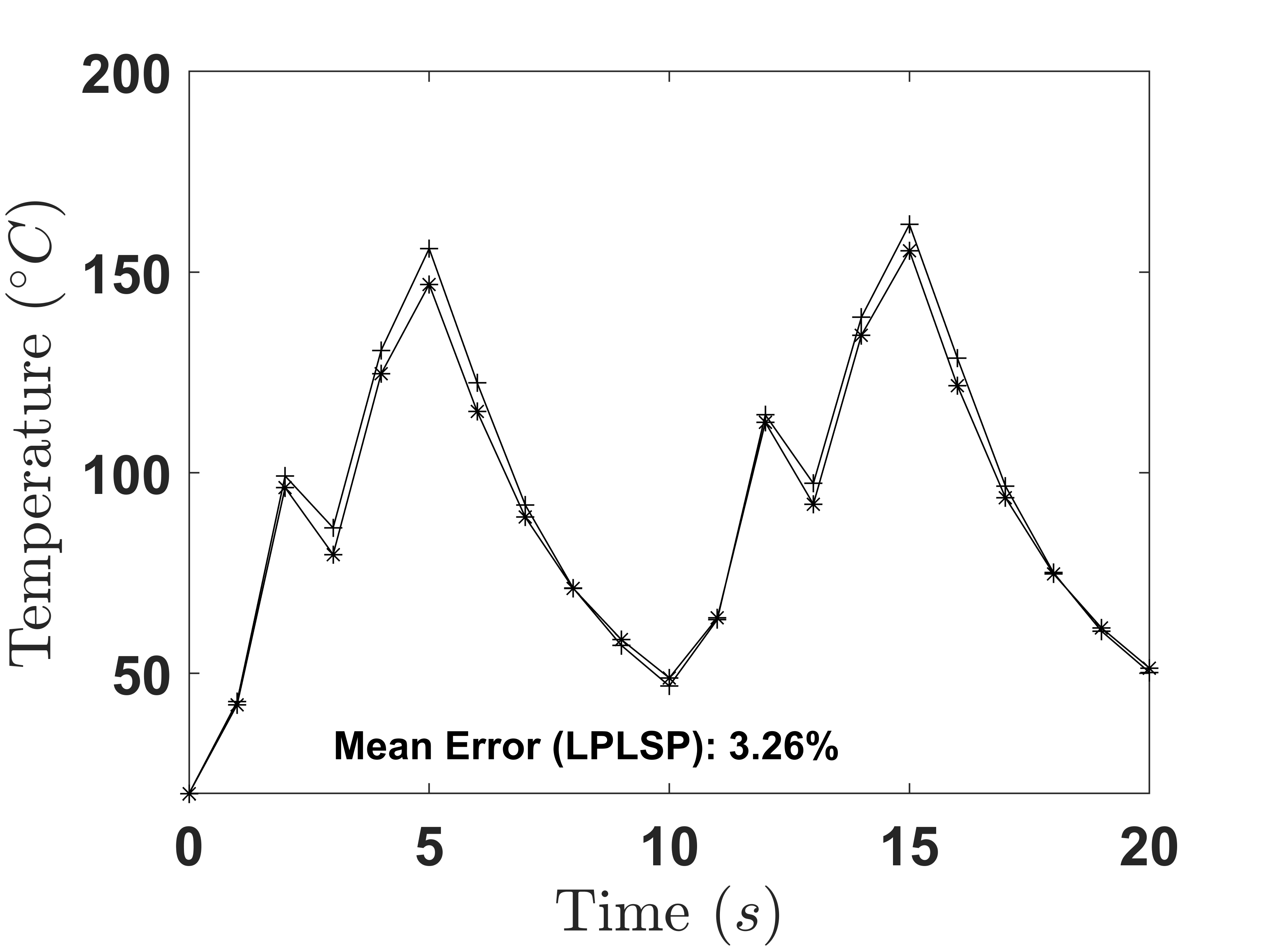}
		\subcaption{MOSFET 3}  
	\end{subfigure} 
	\begin{subfigure}[b]{0.45\textwidth}
		\centering
		\includegraphics[width=\linewidth]{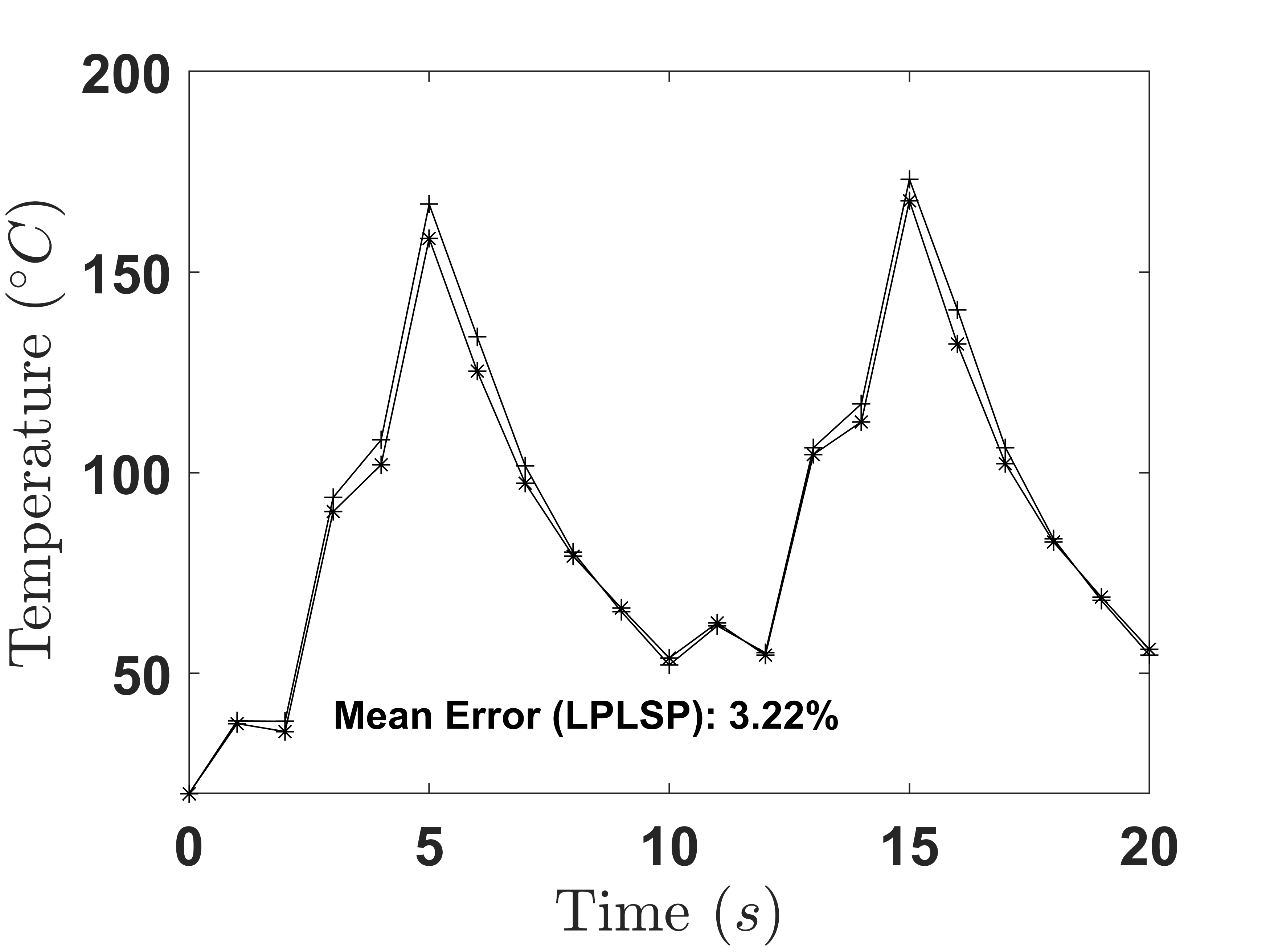}
		\subcaption{MOSFET 4} 
	\end{subfigure} 
	\begin{subfigure}[b]{0.45\textwidth}
		\centering
		\includegraphics[width=\linewidth]{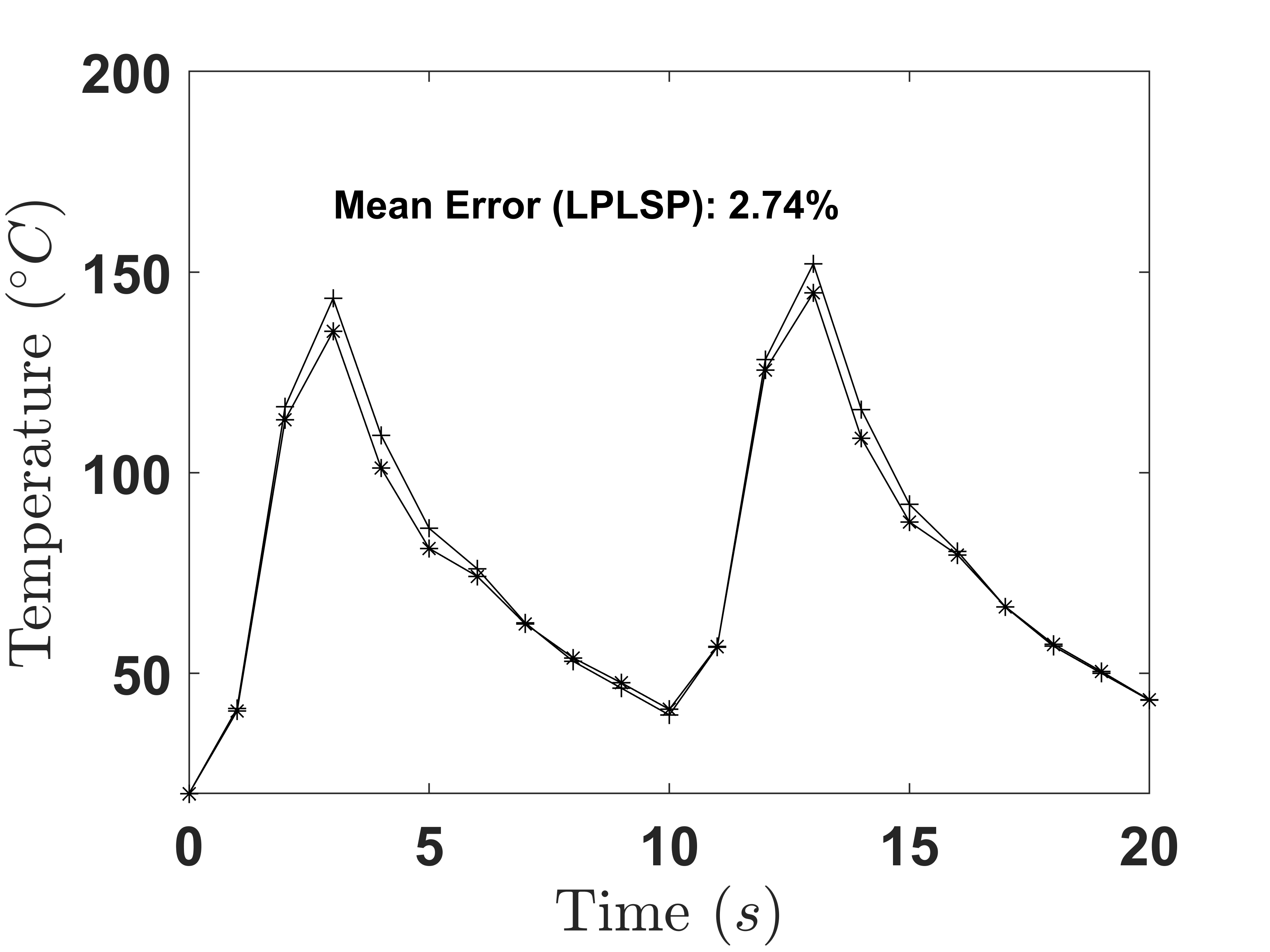}
		\subcaption{MOSFET 5}  
	\end{subfigure}
	\begin{subfigure}[b]{0.45\textwidth}
		\centering
		\includegraphics[width=\linewidth]{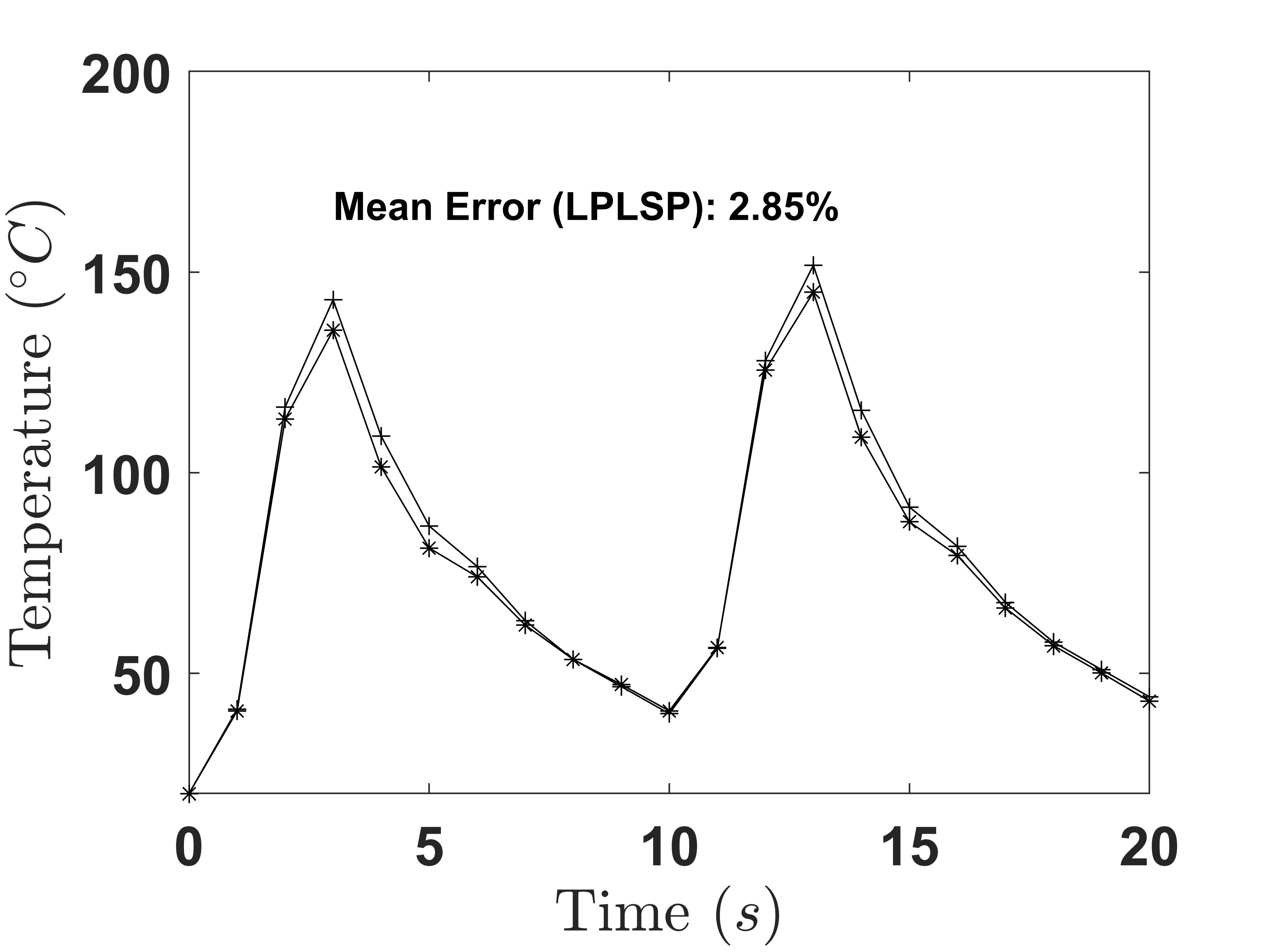}
		\subcaption{MOSFET 6} 
	\end{subfigure} 
	\caption{Comparison of temperatures from simulation $(-*)$ and LPLSP model $(-+)$ for the case of inverter module in a forced convection environment with variable flow velocity profile as presented in Fig.~\ref{fig:velocity}. Mean percentage error $\overline{Error}\% = |\frac{T^{s}_i-T_i}{\max(T^{s}_i)}|$ between the model and simulation data is also presented.}
	\label{fig:FC_uvar}
\end{figure}
\begin{figure}[htbp]\ContinuedFloat
	\centering
	\begin{subfigure}[b]{0.45\textwidth}
		\centering
		\includegraphics[width=\linewidth]{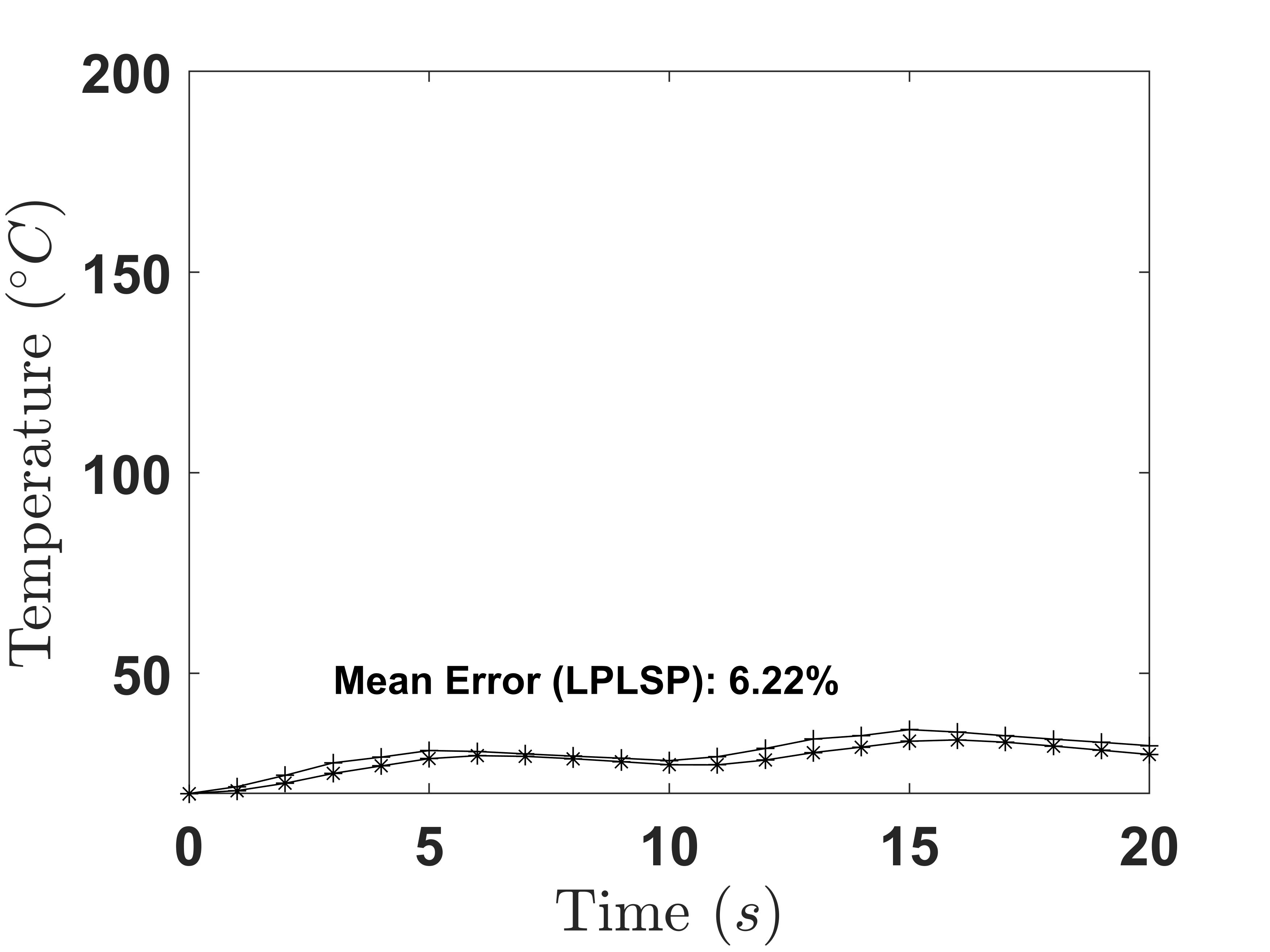}
		\subcaption{PCBA}  
	\end{subfigure} 
	\begin{subfigure}[b]{0.45\textwidth}
		\centering
		\includegraphics[width=\linewidth]{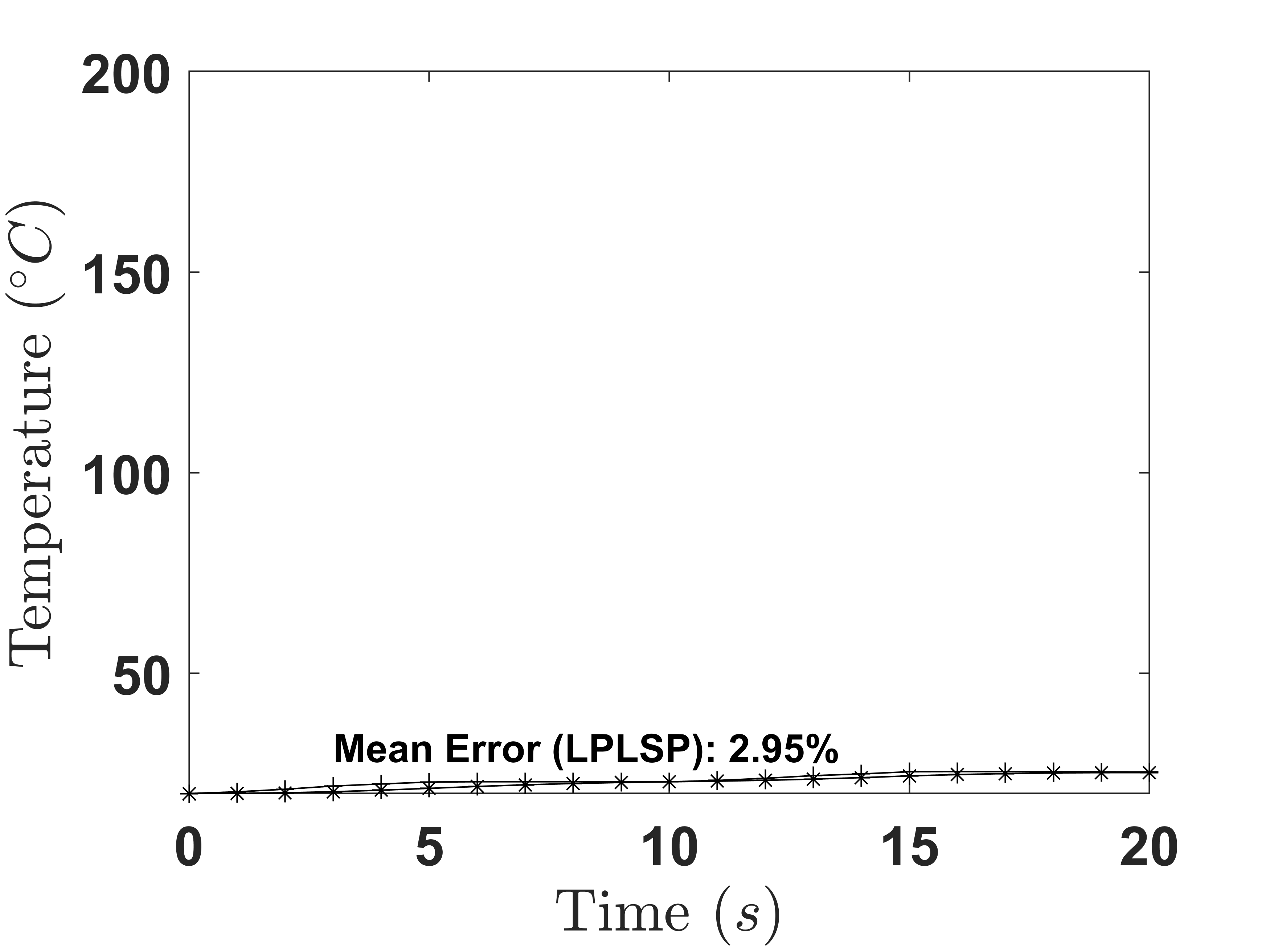}
		\subcaption{Heat Sink} 
	\end{subfigure}
	\caption*{(Figure \ref{fig:FC_uvar} continued)}
\end{figure}
Comparisons between the temperature obtained from CFD simulations, the LPLSP model, and the LTI model are presented for three scenarios: natural convection (Fig.~\ref{fig:NC}), forced convection with constant velocity (Fig.~\ref{fig:FC_const}), and forced convection with variable velocity (Fig.~\ref{fig:FC_uvar}). Since both models are developed based on initial parametric simulations, the computation time is compared in two stages. First, the total time required to perform the parametric simulations is evaluated. Second, the time required to predict the temperature for a new set of inputs is compared. The computation times for both the models and the simulations are summarized in Table~\ref{comptable}. All simulations and model computations were performed on a machine equipped with an Intel Core i7-10850H CPU (12 cores, 2.7 GHz), 32 GB RAM, running Windows 11. All CFD simulations in Ansys Icepak were executed using 6 parallel processors. The LPLSP and LTI model computations were carried out in MATLAB 2019a. LTI models require longer parametric simulation times because they must be run until the system reaches steady-state conditions. While the time taken for a single CFD simulation may be shorter than the combined time for model development and execution, it is important to note that the parametric simulations are performed only once. After the system’s temperature response has been characterized, the resulting model can be reused to simulate any duration and respond to arbitrary changes in input power dissipation or flow velocity.

\begin{table}[h!]
	\centering
	\begin{adjustbox}{max width=\textwidth}
	\begin{tabular}{|c|c|c|c|c|}
		\hline
		Case                                                                                                                                                 & Method     & \begin{tabular}[c]{@{}c@{}}Run time in Ansys Icepak \\ $(s)$\end{tabular} & \begin{tabular}[c]{@{}c@{}}Run time in Matlab \\ $(s)$\end{tabular} & \begin{tabular}[c]{@{}c@{}} Max. error\\ $(\%)$\end{tabular} \\ \hline
		\multirow{3}{*}{\begin{tabular}[c]{@{}c@{}}Natural convection\\ Fig. \ref{fig:NC}\end{tabular}}                                     & Simulation & 900                                                                       & -                                                                   & -                                                                      \\ \cline{2-5} 
		& LPLSP      & 1440                                                                      & 0.563                                                               & 4.74                                                                   \\ \cline{2-5} 
		& LTI        & 4140                                                                      & 6.82                                                                & 25.11                                                                  \\ \hline
		\multirow{3}{*}{\begin{tabular}[c]{@{}c@{}}Forced convection (constant velocity)\\ Fig. \ref{fig:FC_const}\end{tabular}}       & Simulation & 900                                                                       & -                                                                   & -                                                                      \\ \cline{2-5} 
		& LPLSP      & 1440                                                                      & 0.52                                                                & 2.48                                                                   \\ \cline{2-5} 
		& LTI        & 4140                                                                      & 6.02                                                                & 4.97                                                                   \\ \hline
		\multirow{2}{*}{\begin{tabular}[c]{@{}c@{}}Forced convection (variable velocity)\\ Fig. \ref{fig:FC_uvar}\end{tabular}} & Simulation & 730                                                                       & -                                                                   & -                                                                      \\ \cline{2-5} 
		& LPLSP      & 4800                                                                      & 0.338                                                               & 4.22                                                                   \\ \hline
	\end{tabular}
	\end{adjustbox}
	\caption{Comparison of simulation and model run times for the case of inverter module. The run time for the 'Simulation' method in Ansys Icepak represents the total duration required to compute temperatures using full CFD simulation. In contrast, the run times for the LPLSP and LTI methods correspond to the time taken to perform parametric simulations during initial model development. For the LTI method, simulations are executed until steady-state conditions are reached. Although the total run time for the LPLSP method exceeds that of the full simulation, it is important to note that parametric simulations are performed only once. After model development, subsequent evaluations can be executed in MATLAB in under $1 s$. Maximum of the mean absolute error computed for each method is also presented.}
	\label{comptable}
\end{table}

\begin{figure}
	\centering
	\begin{subfigure}[b]{0.32\textwidth}
		\centering
		\includegraphics[width=\linewidth]{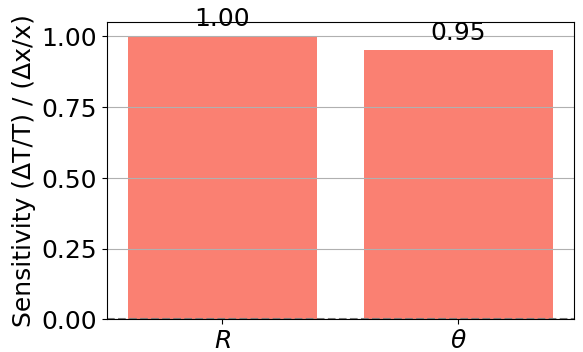}
		\subcaption{For run time $<1 ms$} 
	\end{subfigure}
	\begin{subfigure}[b]{0.32\textwidth}
		\centering
		\includegraphics[width=\linewidth]{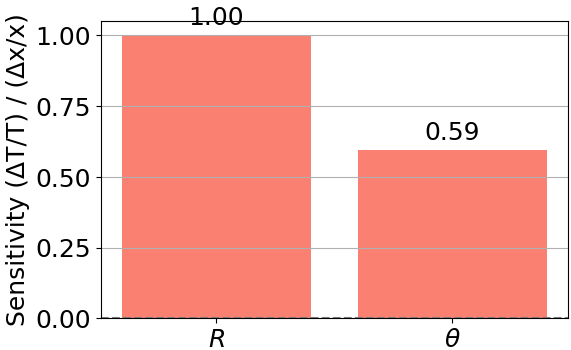}
		\subcaption{For run time $10 ms$}
	\end{subfigure} 
	\begin{subfigure}[b]{0.32\textwidth}
		\centering
		\includegraphics[width=\linewidth]{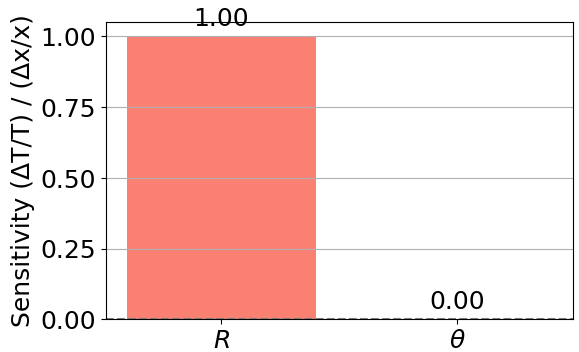}
		\subcaption{For run time $>100 ms$}
	\end{subfigure} \\
	\begin{subfigure}[b]{0.45\textwidth}
		\centering
		\includegraphics[width=\linewidth]{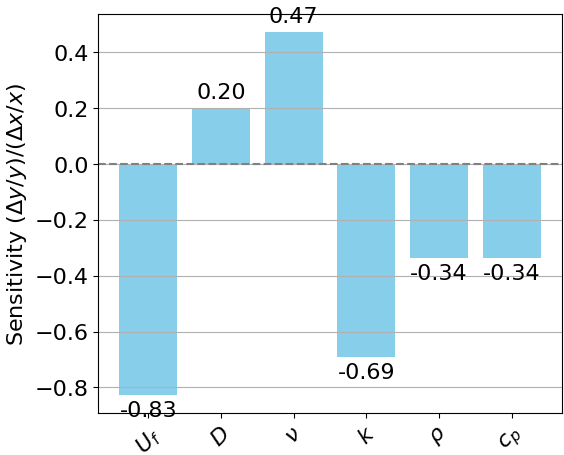}
		\subcaption{Sensitivity of $R$ with respect to physical parameters}
	\end{subfigure}
		\begin{subfigure}[b]{0.45\textwidth}
		\centering
		\includegraphics[width=\linewidth]{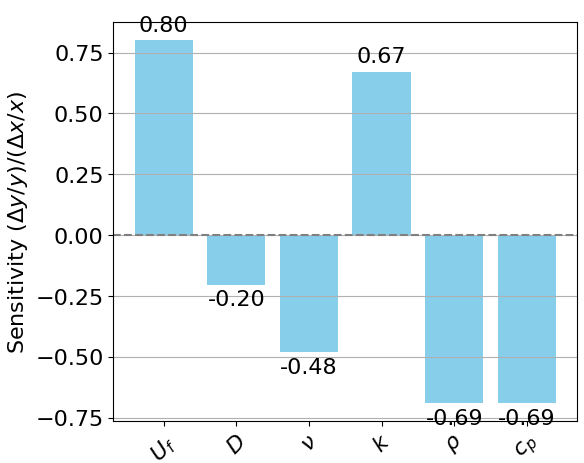}
		\subcaption{Sensitivity of $\theta$ with respect to physical parameters}
	\end{subfigure}
	\caption{Sensitivity of model temperature to parameters $R,C$ and sensitivity of parameters to physical parameters (flow velocity, characteristic length, viscosity, thermal conductivity, density and specific heat) of the system.}
	\label{fig:sensitivity}
\end{figure}

\begin{figure}
	\centering
	\includegraphics[width=\linewidth]{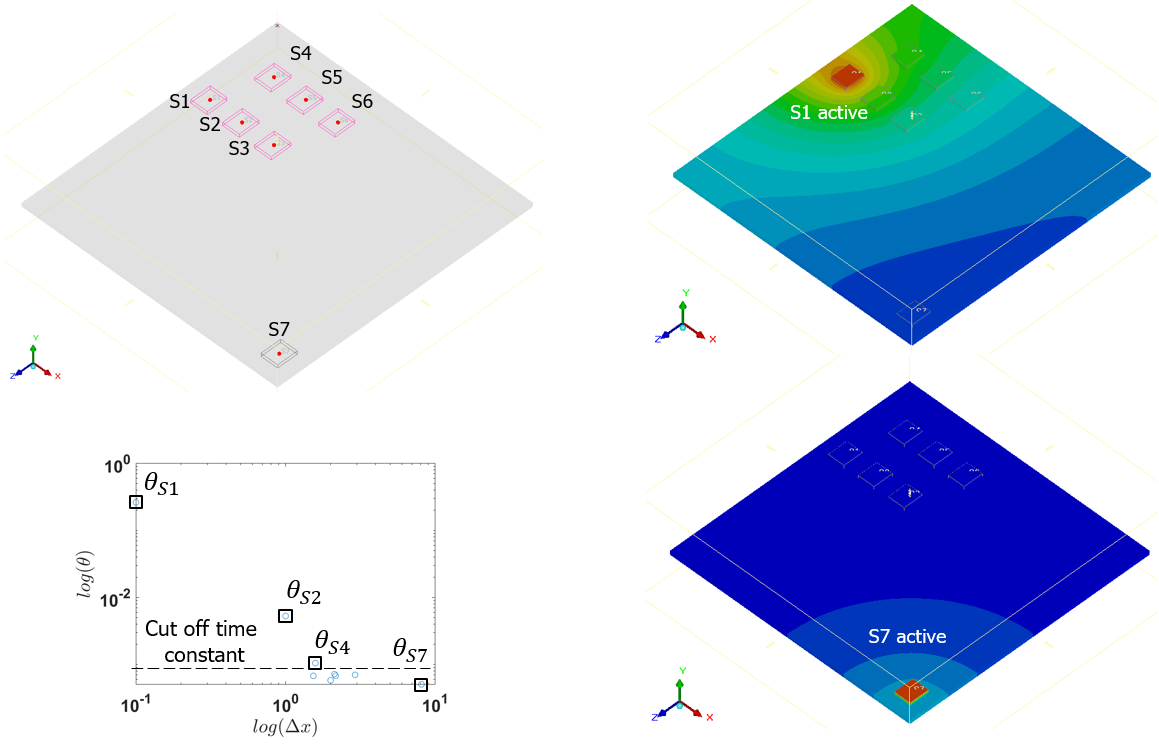}
	\caption{Thermal interaction between sources on the PCB. The influence of S1 on S7 and vice versa is minimal. The time constant measured at S7, with S1 as the heat source, is very small $O(10^{-5})$, indicating negligible thermal coupling. Similarly, the effect of S7 on sources S1–S6 is also insignificant. Therefore, parametric simulations involving S7 can be omitted.}
	\label{fig:kvsx}
\end{figure}

For the case of natural convection, the mean percentage error between the simulation and the LPLSP model is observed to be below $5 \%$, whereas the mean error between the simulation and the LTI model exceeds $20 \%$. This discrepancy arises because of the inherent non-linearity in the physics of natural convection. There is a two way coupling between the flow velocity and the temperature, where the temperature field changes with velocity and velocity is a function of the temperature. Additionally, the flow evolves over time. The fundamental assumption of linearity and time invariance is violated in this case. In this study, the impulse response for the LTI model is obtained by characterizing the steady-state behavior of the temperature. While this steady-state approximation is essential for assuming linear, time-invariant behavior, it becomes problematic in this case. The steady-state temperature curve for natural convection is highly non-linear, and estimating thermal resistance and capacitance from it results in significant modeling error. The LPLSP method avoids this issue by characterizing the temperature response over a shorter time window, where a local linearization assumption holds, leading to improved accuracy. \\
In contrast, for the case of forced convection at a constant flow velocity, the mean percentage error between the simulation and both LTI and LPLSP models are less than $5 \%$. In this case the temperature and the flow velocity are decoupled and the velocity does not evolve dynamically in response to the temperature field. Under these conditions, the convection term in the Navier Stokes and energy equations become linear. Since the flow properties and the driving forces remain constant over time, the resulting flow characteristics also remain constant over time. As a result, the heat transfer behavior remains stable as long as the operating conditions do not change. In this case, the steady-state linearization used by the LTI model aligns well with the local linearization of the LPLSP model, resulting in similar thermal resistance and capacitance values for both models. \\
For the case of forced convection with variable flow velocity, the LTI model is not evaluated since the magnitude of velocity variation is significantly larger than the case of natural convection. Thus the impulse function determined at a fixed velocity will result in inaccurate results. For the LPLSP model the mean percentage error between the temperature from the simulation and the model is still observed to be less than $5 \%$ across all the components. A slightly larger error is observed in the temperature prediction for the PCBA. In this case, the model constants are characterized using a piecewise linear approximation, where $R_{ij} = f(U_f)$ and $\theta_i = f(U_f)$ are defined over small velocity intervals ($U_f = U_0 \pm \delta u$). Specifically, the model constants were computed at discrete flow velocities of $u_x = [1, 5, 10, 20]$ m/s, and linearly interpolated for velocities in between. The overall modeling accuracy can be further improved by refining the velocity resolution, i.e. by computing model constants at shorter intervals which would enable more accurate interpolation and better representation of the system’s velocity-dependent thermal behavior. \\
A sensitivity study is conducted to evaluate the influence of physical parameters of the system on the modeled temperature. The functional relationship in Eq.~\ref{Eq:1BConv} to Eq.~\ref{Eq:h_2} are used as the basis for this study. Geometrical parameters and material properties of the MOSFET and air (as the fluid) are selected as the baseline parameters. A perturbation factor of $\delta = 0.9 (90 \%)$ is applied to each of the parameters. The study is performed in two stages, first the effects of physical parameters such as velocity, characteristic length, viscosity, thermal conductivity, density and specific heat on thermal resistance and time constant are evaluated. Then, the resulting variations in thermal resistance and time constant are used to evaluate their impact on the final temperature. The outcome of the sensitivity study is presented in Fig.~\ref{fig:sensitivity}. It is observed that for short simulation durations on the order of $1 ms$, the final temperature is equally sensitive to both thermal resistance and time constant. For intermediate durations on the order of $10-100 ms$, the temperature becomes more sensitive to thermal resistance. For longer simulation durations exceeding $1 s$, the temperature is predominantly influenced by thermal resistance. Further, for the case study presented in this paper, where model constants in Eq.~\ref{Eq:h_2} are determined from simulation parameters, the thermal resistance is observed to be most sensitive to the flow velocity, followed by thermal conductivity. This indicates that the parameters affecting temperature vary depending on both the application and the time scale of interest. Consequently, the thermal resistances and time constants estimated for modeling are influenced by the duration of the temperature response used for parameter fitting.
While LPLSP method demonstrates high accuracy in this case, it requires a large number of parametric simulations for model development. Specifically, the number of parametric simulations is equal to the number of heat sources in the system. However, results from multiple simulations indicate that components with similar geometry and material properties tend to exhibit similar model parameters when individually powered. For example, the thermal resistance and time constant measured at a component when it is powered ($R_{11}, \theta_{1}$) are comparable to those measured at another similar component when that component is powered ($R_{22}, \theta_{2}$). This can be expressed as $R_{11} \simeq R_{22}$ and $\theta_{1} \simeq \theta_{2}$, and generalized as $R_{kk}, \theta_{k}$ exhibiting similar values across equivalent components. This similarity, however, does not necessarily hold when model parameters are evaluated at a component in response to heat input from a different component (for $R_{ij}, i \neq j$) since the thermal response is dependent on the spatial separation and thermal pathway. Nonetheless, it is observed that components located physically close to each other, but far from the heat source, exhibit similar orders of time constants. Thus, spatial clustering can be exploited by grouping these components and assigning a shared set of time constants. Moreover, when magnitudes of time constants are very low ($O(10^-4)$), the effects of heat source on these components, and vice-versa can be neglected. This can significantly reduce the number of parametric trial runs required to develop the model and further reduce the total computation time. Figure \ref{fig:kvsx} illustrates this concept, where six sources (S1–S6) are located in close proximity on the PCB, while source S7 is located at a farther distance. The thermal time constants of various components are estimated with source S1 active, and their variation with respect to distance from S1 is plotted. The plot clearly demonstrates that the thermal influence of S1 on S7 is significantly weaker compared to its effect on nearby sources such as S2 and S4. Additionally, the influence of S7 on sources S1–S6 is negligible. Consequently, parametric simulations involving S7 can be omitted without loss of accuracy. \\
The LPLSP model, as presented in this work, assumes that the spatial configuration of heat sources and components remains fixed while only the input power dissipation and velocity varies. This assumption is consistent with the motivation behind the method: to enable rapid thermal evaluation for varying electrical and boundary conditions, while keeping the mechanical layout constant. In practical applications such as automotive electronics development, for a particular design iteration, the mechanical design is typically finalized early. Most changes thereafter involve variations in electrical parameters, such as duty cycle, current profiles, and operating conditions, rather than physical relocation of components. In such scenarios, full CFD simulations must be repeated for each new condition, which can be computationally intensive. The LPLSP model addresses this limitation by providing a reduced-order approximation of the CFD simulation, capable of evaluating arbitrary input variations at a fraction of the computational cost. While the model does not explicitly take component location as a parameter, a limited sensitivity study was conducted to evaluate the impact of small positional changes on the model parameters. When the location of a heat-generating component was shifted by approximately $17 \%$ along the $z$-axis, the resulting change in thermal resistance and capacitance was observed to be within $3.12 \%, 3.2 \%$, respectively. Further, the variations in model parameters for the interactions were also limited to a maximum of $8 \%$ This suggests that for small variations in component positions (on the order of $10-20 \%$), the model parameters remain sufficiently robust to provide an accurate approximation. However, it is acknowledged that significant changes in geometry or component layout would necessitate re-evaluation of the model parameters, similar to how full CFD simulations must be re-run in such cases. Thus, the utility of the proposed method lies in its ability to deliver fast and accurate results for a fixed mechanical layout under a wide range of operating conditions, consistent with typical engineering design workflows.
\section{Conclusions}
This paper presents a lumped parameter linear superposition (LPLSP) modeling approach for predicting transient temperature responses in power electronic systems with multiple heat sources and mixed modes of heat transfer. Developed through parametric CFD simulations, the model accurately captures non-linear time-dependent thermal interactions while maintaining a clear, interpretable functional relationship between temperature and relevant physical parameters. Unlike purely data-driven methods, this approach provides a functional relationship between the parameters of the physical system and offers greater insight into the underlying thermal behavior. Comparative analyses against full CFD simulations and LTI-ROM demonstrate that the LPLSP model achieves a mean error below $5 \%$, even in highly nonlinear scenarios such as natural convection and forced convection with variable flow velocities where LTI assumptions break down. The LPLSP method, maintains high fidelity by relying on short-duration transient responses, allowing for effective local linearization even in nonlinear thermal regimes. While the approach requires multiple parametric simulations, strategies to reduce computational overhead were proposed. These include clustering components based on geometric and material similarity, as well as spatial proximity. This allows elimination of terms with small time constant magnitudes, which represent negligible thermal interactions with the source. Nevertheless, some limitations remain. In scenarios where the temperatures of the heat sinks initially exhibit an exponential rise with time constants of order $O(10^{-3})$, followed by a first order exponential rise with time constants of order $O(10^{-1})$, the LPLSP model may result in error due to its purely first order behavior. Further, the model assumes a single phase flow, thus it is not readily applicable to multiphase heat sinks such as vapor chambers or heat pipes. This however could be remedied by incorporating a hybrud sub-model that accounts for latent heat transport and phase change dynamics. This can be achieved by embedding a physics informed surrogate model or a reduced order phase change dynamics within the LPLSP structure. Additionally, the model is primarily empirical in nature and it would benefit from stronger theoretical framework. To reduce the dependence on empirical fitting, the model parameters can be partially derived using analytical correlations and scaling laws based on geometry and material properties. Future work will involve validating the model using experimental data obtained from physical testing. Specifically, by instrumenting a representative power electronic system with thermocouples or IR sensors to measure transient temperature responses under various power and airflow conditions. The experimental setup will be designed to replicate the system and boundary conditions used in CFD simulations, enabling direct comparison with model predictions. Overall, the LPLSP method provides a robust and computationally efficient framework for accurate thermal modeling in complex electronic systems.

\printbibliography

@inproceedings{neel_ecce,
  title={An integrated approach to developing thermal models for automotive electric drives},
  author={Padmanabhan, Neelakantan and Gebregergis, Abraham and Veigas, Santhosh},
  booktitle={2020 IEEE Energy Conversion Congress and Exposition (ECCE)},
  pages={3186--3192},
  year={2020},
  organization={IEEE}
}

@article{zhu2020proper1,
  title={A proper orthogonal decomposition analysis method for transient nonlinear heat conduction problems. Part 1: Basic algorithm},
  author={Zhu, Qiang-Hua and Liang, Yu and Gao, Xiao-Wei},
  journal={Numerical Heat Transfer, Part B: Fundamentals},
  volume={77},
  number={2},
  pages={87--115},
  year={2020},
  publisher={Taylor \& Francis}
}

@article{zhu2020proper2,
  title={A proper orthogonal decomposition analysis method for transient nonlinear heat conduction problems. Part 2: Advanced algorithm},
  author={Zhu, Qiang-Hua and Liang, Yu and Gao, Xiao-Wei},
  journal={Numerical Heat Transfer, Part B: Fundamentals},
  volume={77},
  number={2},
  pages={116--137},
  year={2020},
  publisher={Taylor \& Francis}
}

@article{neel_pod,
  title={Reduced order model of a convection-diffusion equation using Proper Orthogonal Decomposition},
  author={Padmanabhan, Neelakantan},
  journal={arXiv preprint arXiv:2303.07176},
  year={2023}
}

@article{neel_chapter,
  title={Recent Advancements in Large Eddy Simulations of Compressible Real Gas Flows},
  author={Padmanabhan, Neelakantan},
  year={2024},
  publisher={IntechOpen}
}

@article{pod_theory,
  title={The proper orthogonal decomposition in the analysis of turbulent flows},
  author={Berkooz, Gal and Holmes, Philip and Lumley, John L},
  journal={Annual review of fluid mechanics},
  volume={25},
  number={1},
  pages={539--575},
  year={1993},
  publisher={Annual Reviews 4139 El Camino Way, PO Box 10139, Palo Alto, CA 94303-0139, USA}
}

@article{haider2008proper,
  title={A proper orthogonal decomposition based system-level thermal modeling methodology for shipboard power electronics cabinets},
  author={Haider, Syed I and Burton, Ludovic and Joshi, Yogendra},
  journal={Heat transfer engineering},
  volume={29},
  number={2},
  pages={198--215},
  year={2008},
  publisher={Taylor \& Francis}
}

@InProceedings{CTN1,
author="Denny, Allen
and Kirloskar, Neelkanth
and Ponangi, Babu Rao
and Joseph, Rex
and Krishna, V.",
editor="Krishna, V.
and Seetharamu, K. N.
and Joshi, Yogendra Kumar",
title="Electro-Thermal Model for Field Effective Transistors",
booktitle="Recent Advances in Hybrid and Electric Automotive Technologies",
year="2022",
publisher="Springer Nature Singapore",
address="Singapore",
pages="277--284",
isbn="978-981-19-2091-2"
}

@article{CTN2,
  title={A high-precision adaptive thermal network model for monitoring of temperature variations in insulated gate bipolar transistor (IGBT) modules},
  author={An, Ning and Du, Mingxing and Hu, Zhen and Wei, Kexin},
  journal={Energies},
  volume={11},
  number={3},
  pages={595},
  year={2018},
  publisher={MDPI}
}

@article{CTN3,
  title={Power Semiconductor Junction Temperature and Lifetime Estimations: A Review},
  author={Morel, Cristina and Morel, Jean-Yves},
  journal={Energies},
  volume={17},
  number={18},
  pages={4589},
  year={2024},
  publisher={MDPI}
}

@article{CTN4,
title = {Behavioral electrothermal modeling of MOSFET for energy conversion circuits simulation using MATLAB/Simulink},
journal = {Microelectronics Reliability},
volume = {154},
pages = {115340},
year = {2024},
issn = {0026-2714},
doi = {https://doi.org/10.1016/j.microrel.2024.115340},
url = {https://www.sciencedirect.com/science/article/pii/S0026271424000209},
author = {Mohamed Baghdadi and Elmostafa Elwarraki and Imane {Ait Ayad} and Naoual Mijlad},
}

@article{fostercauer,
  title={A new approach to the dynamic thermal modelling of semiconductor packages},
  author={Masana, FN},
  journal={Microelectronics Reliability},
  volume={41},
  number={6},
  pages={901--912},
  year={2001},
  publisher={Elsevier}
}

@article{krylov1,
  title={Theory of Krylov subspace methods based on the Arnoldi process with inexact inner products},
  author={Su, Meng and Wen, Chun and Shen, Zhao-Li and Serra-Capizzano, Stefano},
  journal={NHM},
  volume={20},
  number={1},
  pages={15--34},
  year={2025}
}

@article{krylov2,
  title={KRYLOV SUBSPACE MODEL ORDER REDUCTION OF LARGE SCALE FINITE ELEMENT DYNAMICAL SYSTEMS},
  author={Sindler, Jaroslav and Sulitka, Matej},
  journal={MM Science Journal},
  year={2013}
}

@article{krylov3,
  title={The extended and asymmetric extended krylov subspace in moment-matching-based order reduction of large circuit models},
  author={Stoikos, Pavlos and Garyfallou, Dimitrios and Floros, George and Evmorfopoulos, Nestor and Stamoulis, George},
  journal={arXiv preprint arXiv:2204.02467},
  year={2022}
}

@article{krylov4,
  title={Practical issues of model order reduction with Krylov-subspace methods},
  author={Heres, PJ and Schilders, WHA},
  year={2004},
  publisher={Technische Universiteit Eindhoven}
}

@phdthesis{krylov5,
  title={Krylov subspace techniques for model reduction and the solution of linear matrix equations},
  author={Ahmad, Mian Ilyas},
  year={2011},
  school={Imperial College London}
}

@article{BT1,
  title={Balanced truncation model reduction for linear time-varying systems},
  author={Lang, Norman and Saak, Jens and Stykel, Tatjana},
  journal={Mathematical and computer modelling of dynamical systems},
  volume={22},
  number={4},
  pages={267--281},
  year={2016},
  publisher={Taylor \& Francis}
}

@article{BT2,
  title={Balanced truncation model reduction of second-order systems},
  author={Reis, Timo and Stykel, Tatjana},
  journal={Mathematical and Computer Modelling of Dynamical Systems},
  volume={14},
  number={5},
  pages={391--406},
  year={2008},
  publisher={Taylor \& Francis}
}

@article{ML1,
  title={Review of machine learning techniques for power electronics control and optimization},
  author={Bahrami, Maryam and Khashroum, Zeyad},
  journal={arXiv preprint arXiv:2310.04699},
  year={2023}
}

@article{ML2,
  title={Artificial intelligence for power electronics in electric vehicles: challenges and opportunities},
  author={Paret, Paul and Finegan, Donal and Narumanchi, Sreekant},
  journal={Journal of Electronic Packaging},
  volume={145},
  number={3},
  pages={034501},
  year={2023},
  publisher={American Society of Mechanical Engineers}
}

@misc{ML3,
  title = {Challenges in Enhancing Power Electronic Systems with Artificial Intelligence},
  howpublished = {\url{https://www.utmel.com/blog/categories/technology/challenges-in-enhancing-power-electronic-systems-with-artificial-intelligence}},
  note = {Accessed: 2025-06-26}
}

@inproceedings{neel_ieee,
  title={A Transient Thermal Model for Power Electronics Systems},
  author={Padmanabhan, Neelakantan},
  booktitle={SoutheastCon 2024},
  pages={1294--1299},
  year={2024},
  organization={IEEE}
}

@misc{RC_tr1,
  title = {Thermal Resistance Theory and Practice},
  author = {Infineon Technologies},
  year = {2000}
}

@article{RC_tr2,
  title={Modeling the transient response of thermal circuits},
  author={Silva, Daniel},
  journal={Applied Sciences},
  volume={12},
  number={24},
  pages={12555},
  year={2022},
  publisher={MDPI}
}

@misc{neel_cfd2,
	title={On High Pressure Real Gas Turbulent Mixing Jets, 2017. All Dissertations. 1981},
	author={Padmanabhan, N},
	year={1981}
}



\end{document}